\RequirePackage{fix-cm} 
\documentclass[a4paper,twoside,12pt,reqno]{amsart}

\usepackage{fixltx2e}     

\usepackage[english]{babel}
\usepackage[latin1]{inputenc}

\usepackage{indentfirst,verbatim}

\usepackage{amsmath,amsfonts,amssymb,amsgen,amsbsy,eucal,mathrsfs,dsfont}
\usepackage{stmaryrd}

\usepackage[square,numbers]{natbib}
\usepackage{a4wide,verbatim}

\usepackage{epsfig,rotating,color}
\usepackage{multicol}
\usepackage{tikz,pgfplots}
\usepackage{graphicx}
\usepackage{mathtools}

\usepackage{caption}

\usepackage{hyperref}

\usepackage[square,numbers]{natbib}
\hypersetup{
colorlinks=true, 
breaklinks=true, 
urlcolor= blue, 
linkcolor= blue, 
bookmarksopen=true, 
pdftitle={}, 
pdfauthor={}, 
pdfsubject={} 
}

\newcommand{\ud}{\mathrm{d}}
\newcommand{\uD}{\mathrm{D}}

\newcommand{\phim}{{\phi}}

\newcommand{\half}{{\textstyle{1\over2}}}

\newtheorem{thm}{Theorem}
\newtheorem{lem}{Lemma}

\newtheorem{rem}{Remark}

\usepackage{enumitem}
\newlist{steps}{enumerate}{1}
\setlist[steps, 1]{label = Step \arabic*:}

\newcommand{\eqdef}{\stackrel{\text{\tiny{def}}}{=}}

\title[Regularised barotropic Euler equation]{\bf Hamiltonian regularisation of the 
unidimensional barotropic Euler equations}

\author[GUELMAME et al.]{Billel Guelmame, Didier Clamond and St\'ephane Junca}

\newcommand{\nfont}{\fontshape{n}\selectfont}

\address{({\nfont\textbf{Billel Guelmame}})  LJAD,  Inria \& CNRS,  Universit\'e C\^ote d'Azur, France.} 
\email{billel.guelmame@univ-cotedazur.fr}

\address{({\nfont\textbf{Didier Clamond}})  LJAD,  CNRS,  Universit\'e C\^ote d'Azur,  France.} 
\email{didier.clamond@univ-cotedazur.fr}

\address{({\nfont\textbf{St\'ephane Junca}})  LJAD,  Inria \& CNRS,   Universit\'e C\^ote d'Azur, France.} 
\email{stephane.junca@univ-cotedazur.fr}

\let\oldtocsection=\tocsection
 
\let\oldtocsubsection=\tocsubsection

\renewcommand{\tocsection}[2]{\hspace{0em}\oldtocsection {#1}{#2}}
\renewcommand{\tocsubsection}[2]{\hspace{2em}\oldtocsubsection{#1}{#2}}

\begin{document}

\maketitle

\begin{abstract}
Recently, a Hamiltonian regularised shallow water (Saint-Venant) system has been introduced by \citet{ClamondDutykh2018a}. This system is Galilean invariant, linearly non-dispersive and conserves 
formally an $H^1$-like energy. In this paper, we extend this regularisation in two directions. 
First, we consider the more general barotropic Euler system, the shallow water equations being 
formally a very special case. Second, we introduce a class regularisations, showing thus that 
this Hamiltonian regularisation of \citep{ClamondDutykh2018a} is not unique. 
Considering the high-frequency approximation of this regularisation, we obtain a new two-component Hunter--Saxton system.
We prove that both systems --- the regularised barotropic Euler system  and the two-component Hunter--Saxton system --- are locally (in time) well-posed, and, if singularities appear in finite time, they are necessary in the first derivatives.
\end{abstract}

\medskip


\medskip


\tableofcontents

\section{Introduction}

The barotropic Euler system is a quasilinear system of partial differential equations that can be 
used to describe many phenomenae in fluid mechanics. Denoting the time and the spacial coordinate 
by the independent variables $t$ and $x$, respectively, and denoting the density, the velocity and 
the pressure by the dependent variables $\rho(t,x)>0$, $u(t,x)$ and $P(\rho)$, respectively,
the conservation of mass and momentum yield
\begin{equation}\label{bE}
\rho_t\ +\,\left[\,\rho\,u\,\right]_x\ =\ 0, \qquad
u_t\ +\ u\,u_x\ +\ P_x/\rho\ =\ 0,
\end{equation}
where subscripts denote partial derivatives.  The system \eqref{bE} is the unidimensional Euler 
equations for compressible perfect fluids. With $P=P(\rho)$ we have the barotropic Euler equations. 
Many equations of state for compressible fluids can be found in the literature \citep{WagnerPruss2002}, 
the most commonly used being probably the isentropic motion (Appendix \ref{appisent}). 
If the pressure $P$ is an increasing function of $\rho$, the system \eqref{bE} is hyperbolic 
and, even for smooth initial data, the solutions may develop shocks in finite time. 

In order to avoid the formation of shocks, several regularisations have been proposed, for example 
by adding a ``small'' artificial viscosity and/or dispersive terms \citep{bianchini2005vanishing,
crandall1983viscosity,HayesLeFloch2000,kondo2002zero,LaxLevermore1983-1,LaxLevermore1983-2,
LaxLevermore1983-3,vonNeumannRichtmyer1950}. 
The artificial viscosity leads to a loss of the energy everywhere and the dispersive terms lead to 
high-frequency oscillations which may cause problems in numerical computations. Other regularisations 
of Leray-type (for inviscid Burgers equation, isentropic Euler system and others) have also been proposed and 
studied in \citep{BhatFetecau2006,BhatFetecau2009a,BhatFetecau2009b,BhatEtAl2007,CamassaEtAl2010,
NorgardMohseni2010a,NorgardMohseni2010b}. 
Those regularisations do not conserve the energy and the limit solution fails (in general) to satisfy 
the Lax entropy condition \cite{BhatFetecau2009a}. 

For the Saint-Venant (shallow water) system of hyperbolic equations, \citet{ClamondDutykh2018a} 
proposed the dispersionless non-diffusive regularised  system
\begin{subequations}\label{rSV}
\begin{align}
h_t\ +\,\left[\, h\,u\,\right]_x\ =\ 0, \\
\left[\,h\,u\,\right]_t\ +\,\left[\,h\,u^2\,+\, \half\, g\, h^2\,
+\,\epsilon\,\mathscr{R}\,\right]_x\ =\ 0,\\
\mathscr{R}\, \eqdef\, 2\, h^3\, u_x^2\, -\, h^3 [u_t\, +\, u\, u_x\, +\ g\, h_x]_x\, -\, 
\half\, g\, h^2\, h_x^2,
\end{align}
\end{subequations}
where $h$ is the total water depth of the fluid, $g>0$ is the gravitational acceleration and 
$\varepsilon>0$ is a dimensionless parameter. If $\mathscr{R}=0$ the classical shallow water 
equations are recovered. The suitable regularising term $\mathscr{R}$ has been obtained 
applying a Leray-like regularisation to the variational principle leading to the shallow water 
equations; it is therefore called {\em Hamiltonian regularisation} since the equations \eqref{rSV} 
can be rewritten in a non-canonical Hamiltonian form (see \citep{ClamondDutykh2018a} and 
\citep{PuEtAl2018} for details). This regularisation conserves an $H^1$-like energy for 
smooth solutions and it has the same velocity of propagation of singularities as the classical Saint-Venant system. Weak 
singular shocks of \eqref{rSV} have been studied by \citet{PuEtAl2018}. Also, local (in time) 
well-posedness and existence of blowing-up solutions using Ricatti-type equations have been 
proved in \citep{liu2019well}. The global well-posedness and a mathematical study of the case 
$\varepsilon \to 0$ remain open problems. 

Recently, inspired by \citep{ClamondDutykh2018a} and with the same properties as \eqref{rSV}, 
a similar regularisation has been proposed for the inviscid Burgers equation in \citep{guelmame2020global} 
and for general scalar conservation laws in \citep{guelmame2020hamiltonian}, where solutions 
exist globally (in time) in $H^1$, those solutions converging to solutions of the classical 
equation when $\varepsilon\to 0$ at least for a short time \citep{guelmame2020hamiltonian,
guelmame2020global}. The regularised Saint-Venant system \eqref{rSV} has been also generalised 
for shallow water equations with uneven bottoms \cite{clamond2019hamiltonian}.

The classical Saint-Venant system (letting $\varepsilon\to0$ in \eqref{rSV}) is, formally, a special 
case of the barotropic Euler system \eqref{bE} obtained substituting the density $\rho$ for the total 
depth $h$ (i.e., $h\mapsto\rho$) and substituting the pressure $P$ for $\half g h^2$. One thus 
obtains the special isentropic Euler equation with pressure $P\propto \rho^2$ (i.e., 
classical isentropic Euler equation with $\gamma=2$, see Appendix \ref{appisent}). 
One goal of this paper is to extend the regularised system \eqref{rSV} to the barotropic Euler 
system \eqref{bE} with $P=P(\rho)$, in particular for $P\propto\rho^\gamma$ for any $\gamma>0$. 

The regularised Saint-Venant equations \eqref{rSV} were originally derived in \citep{ClamondDutykh2018a} 
from a modified version of the Serre--Green--Naghdi equations \citep{ClamondEtAl2017}. We show 
here that it is actually only a special case of this type of regularisation. A second 
goal of this paper is then to introduce a class of regularisations applied to the general 
barotropic Euler equations \eqref{bE} as in 
\citep{ClamondDutykh2018a}, that is preserving the same properties. In Section \ref{secregbar} below, 
modifying the Lagrangian of \eqref{bE}, we obtain the regularised barotropic Euler (rbE) system
\begin{subequations}\label{rE}
\begin{align} 
\rho_t\ +\,\left[\,\rho\,u\,\right]_x\ =\ 0, \label{rEa} \\
\left[\,\rho\,u\,\right]_t\ +\,\left[\,\rho\,u^2\,+\,\rho\,\mathscr{V}'\,-\,\mathscr{V}\,
+\,\epsilon\,\mathscr{R}\,\right]_x\ =\ 0, \label{rEb}\\
\mathscr{R}\ \eqdef\,\left(\/\rho^2\/\mathscr{A}'\/\right)'u_x^{\,2}\ -\ 2\,\rho\,\mathscr{A}'
\left[\,u_t\,+\,u\,u_x\,+\,\varpi_x\,\right]_x\ +\ \,\left(\/\rho\/\mathscr{V}''\/
/\/\mathscr{A}'\right)'\,\mathscr{A}_x^{\,2}, \label{rEc}
\end{align}
\end{subequations}
where primes denote derivatives with respect to $\rho$, $\mathscr{V}''(\rho)=P'(\rho)/\rho$ and 
$\mathscr{A}$ is a smooth increasing function of $\rho$. We show in this paper that the system 
\eqref{rE} is non-dispersive, non-diffusive, it has a Hamiltonian structure, it has the same velocity of propagation of singularities as \eqref{bE} and, for all $\mathscr{A}$ smooth increasing function of $\rho$, smooth solutions 
of \eqref{rE} conserve an $H^1$-like energy.
A study of steady solutions of rbE is also done, which covers the traveling waves due to the 
Galilean invariance of rbE.

Introducing the linear Sturm--Liouville operator 
\begin{equation}\label{op:SL}
\mathcal{L}_\rho\ \eqdef\ \rho\ -\ 2\, \epsilon\,
\partial_x\, \rho\, \mathscr{A}'\,\partial_x,
\end{equation} 
and applying $\mathcal{L}_\rho^{-1}$ on \eqref{rEb} 
(the invertibility of $\mathcal{L}_\rho$ is insured by Lemma \ref{Inverseesitimates} below), the 
system \eqref{rE} becomes
\begin{subequations}\label{rE2}
\begin{align} 
\rho_t\ +\,\left[\,\rho\,u\,\right]_x\ &=\ 0, \\
u_t\ +\ u\,u_x\ +\ P_x/\rho\ &=\ -\,\epsilon\,\mathcal{L}^{-1}_\rho\partial_x \left\{\,\left(\/\rho^2\/\mathscr{A}'\/\right)'u_x^{\,2}\, 
+\left(\/\rho\/\mathscr{V}''\//\/\mathscr{A}'\right)'\,{\mathscr{A}'}^{\/2}\, \rho_x^{\,2}\,\right\}.
\end{align}
\end{subequations}
The reason for applying $\mathcal{L}_\rho^{-1}$ is to remove the derivative with respect to $t$ and the 
high-order derivatives with respect to $x$ appearing in \eqref{rEc}. 
The form \eqref{rE2} is then more convenient to obtain the local well-posedness of rbE. 
Following \citep{Alinhac,
AlinhacGerard,
Israwi2011,
liu2019well,
majda2012compressible}, we prove that 
if the initial data is an $H^s$ perturbation of a constant state (with $s\geqslant 2$, and $\rho \geqslant 
\rho^{*}>0$), then \eqref{rE2} is locally well-posed. 
The same proof is used to prove the local existence of periodic solutions of the generalised two-component Hunter--Saxton 
system
\begin{subequations}\label{HS}
\begin{align}
\rho_t\ +\,\left[\,\rho\,u\,\right]_x\ &=\ 0, \\ 
u_t\,+\,u\,u_x\,+\,P_x/\rho\ 
&=\, \partial_x^{-1}\,  \left\{ \left(\,1\,+\,\frac{\rho\,\mathscr{A}''}{2\,\mathscr{A}'}
\,\right)u_x^{\,2}\ +\,\left(\,\frac{\left(\rho\/\mathscr{V}''\right)'}{2\,\rho}\,-\,\frac{\mathscr{V}''\,
\mathscr{A}''}{2\,\mathscr{A}'}\,\right)\rho_x^{\,2} \right\},
\end{align}
\end{subequations}
that can be obtained by formally taking $\varepsilon\to \infty$ in \eqref{rE2}.
We show also in Section \ref{secsr} that if $\mathscr{A}(\rho) = -A/\rho$ where $A>0$, then the generalised Hunter--Saxton 
system \eqref{HS} is formally equivalent to the so-called variational wave equation 
\begin{equation*}
v_{tt}\ -\ c(v)\, \left( c(v)\, v_x \right)_x\ =\ 0.
\end{equation*}

This paper is divided in two parts. 
A first part (Sections \ref{secbE} \& \ref{secregbar}) presents the physical motivations of the 
regularised barotropic Euler system and its properties.
A second part (Sections \ref{secapp} \& \ref{secHS}) consists of mathematical proofs for  existence 
results. 
Shortly speaking, the first part is more physical oriented, while the second one is more mathematical.    
More specifically, the content of the paper is organised as follows.
In Section \ref{secbE}, we recall some classical properties of the 
barotropic Euler system \eqref{bE}. Section \ref{secregbar} is devoted to derive the regularised system 
\eqref{rE}, to study its properties and steady motions. 
In Section \ref{secapp}, we prove the local well-posedness and a blow-up criteria of \eqref{rE2}.
In Section \ref{secHS}, the generalised Hunter--Saxton system \eqref{HS} is introduced, and a well-posedness theorem is given.
A special choice of the regularising function $\mathscr{A}$, with some 
interesting properties, is briefly discussed in Section \ref{secsr}.

\section{Equations for barotropic perfect fluids}\label{secbE}

Let us recall the conservation of mass and momentum for perfect fluids in Eulerian description of motion
\begin{subequations}\label{eqEul}
\begin{align}\label{eqEula}
\rho_t\ +\,\left[\,\rho\,u\,\right]_x\ &=\ 0, \\ \label{eqEulb}
u_t\ +\ u\,u_x\ +\ P_x/\rho\ &=\ 0.
\end{align}
\end{subequations}

Note that (conservative) body forces, if present, are incorporated into the definition of the pressure $P$.
In the special case of barotropic motions \citep{Serrin1959} --- i.e., when $\rho=\rho(P)$ or $P=P(\rho)$ 
--- it is convenient to introduce the so-called {\em specific enthalpy} \citep[\S3.3]{Dolzhansky2013} 
$\varpi$ such that
\begin{align} \label{defvarpi}
\varpi\ =\ \int\frac{\ud\/P}{\rho(P)}\ =\ \int\frac{\ud\,P(\rho)}{\ud\/\rho}\,
\frac{\ud\/\rho}{\rho} \qquad\Rightarrow\qquad
\rho\,\ud\varpi\ =\ \ud P, \qquad \partial_x \varpi\ =\ \frac{\partial_x P}{\rho}. 
\end{align}
$\varpi$ being an antiderivative of $1/\rho(P)$, it is defined modulo an additional arbitrary 
integration constant, so the value of $\varpi$ can be freely chosen on a given isobaric surface 
$P=\text{constant}$ (thus providing a gauge condition for the specific enthalpy). 
The relation (\ref{defvarpi}) can also be written in the reciprocal form $P(\varpi)=
\int\rho(\varpi)\,\ud\/\varpi$, thence $P$ is a known function of $\varpi$ and, obviously, 
$P=\rho\/\varpi$ if the density is constant. The speed of sound $c_s$ is defined by
\begin{equation*}\label{defcs}
c_s\ \eqdef\ \left[\,{\ud\/\rho}\//\/{\ud\/P}\,\right]^{-{1\over2}}\ =\ 
\left[\,\rho^{-1}\,{\ud\/\rho}\//\/{\ud\/\varpi}\,\right]^{-{1\over2}}.
\end{equation*} 
From this definition we have $\rho\,\ud\varpi=c_s^{\,2}\,\ud\rho$, thence
with the mass conservation (\ref{eqEula})
\begin{equation}\label{eqvpiucs}
\uD_t\,\varpi\ =\ \rho^{-1}\,c_s^{\,2}\,\uD_t\,\rho\ =\ -\/c_s^{\,2}\ u_x,
\end{equation}
where $\uD_t\eqdef\partial_t+u\partial_x$ is the temporal derivative following the motion.
The relation (\ref{eqvpiucs}) is of special interest when $c_s$ is constant.
Many equations of state for compressible fluids can be found in the literature \citep{WagnerPruss2002}. 
Isentropic motions are of special interest so their equation of state is given in Appendix \ref{appisent}.

\subsection{Cauchy--Lagrange equation}

For barotropic fluids, the momentum equation (\ref{eqEulb}) becomes
\begin{equation}\label{mompi}
u_t\ +\ u\,u_x\ +\ \varpi_x\ =\ 0,
\end{equation}
and introducing a velocity potential $\phi$ such that $u=\phi_x$, the equation (\ref{mompi}) 
is integrated into a Cauchy--Lagrange equation
\begin{equation}\label{eqbern}
\phi_t\ +\ \half\,\phi_x^{\,2}\ +\ \varpi\ =\ K(t)\ \equiv\ 0,
\end{equation}
where $K(t)$ is an integration constant that can be set to zero without loose of generality 
(gauge condition for the velocity potential).

\subsection{Conservation laws}

For regular solutions, secondary conservation laws can be derived from (\ref{eqEul}), e.g.,
\begin{align}
%
\left[\/\rho\/u\/\right]_t\,+\left[\/\rho\/u^2\/+\/\rho\/\mathscr{V}'\/-\/\mathscr{V}\/\right]_x\, &=\, 0,\\
\left[\/u\/\right]_t\,+\left[\/\half\/u^2\/+\/\mathscr{V}'\/\right]_x\, &=\, 0, \\
\left[\/\half\/\rho\/u^2\/+\/\mathscr{V}\/\right]_t\, +\left[\left(\/\half\/\rho\/u^2\/+\/
\rho\/\mathscr{V}'\/\right)u\/\right]_x\, &=\, 0,
%
%
\end{align}
where
\begin{gather*}
\mathscr{V}\, \eqdef\, \int\varpi\,\ud\/\rho.
\end{gather*}
Actually, since the barotropic Euler system is a $2 \times 2$ strictly hyperbolic system, an infinite number of conservation laws can be derived \citep{Dafermos2016,Laxbook}. 

\subsection{Jump conditions}
 
The Euler equations admit weak solutions. For discontinuous $\rho$ and $u$, the 
Rankine--Hugoniot conditions for the mass and momentum conservation are 
\begin{equation}\label{jumpeuler}
\left(u\,-\,\dot{s}\right)\llbracket\,\rho\,\rrbracket\,+\,\rho\,\llbracket\,u\,\rrbracket\, =\ 0,
\qquad
\left(u\,-\,\dot{s}\right)\llbracket\,u\,\rrbracket\,+\,\varpi'\,\llbracket\,\rho\,\rrbracket\, =\ 0,
\end{equation}  
where $\varpi'=\ud\varpi/\ud\rho$, $\dot{s} \eqdef \ud s/\ud t$ is the velocity of the propagation of the weak discontinuity located 
at $x=s(t)$ and $\llbracket\/f\/\rrbracket\eqdef f(x\!=\!s^+)-f(x\!=\!s^-)$ denotes the jump across 
the discontinuity for any function $f$. The Rankine--Hugoniot conditions (\ref{jumpeuler}) yield at once 
the velocity 
\begin{equation}\label{velocityeuler}
\dot{s}(t)\ =\ u\ \pm\ \sqrt{\,\rho\,\varpi'\,} \qquad \text{ at }\quad x\,=\,s(t).
\end{equation}
A goal of the present work is to derive a regularisation of the Euler equation that 
preserves exactly the same velocity of the propagation of the weak discontinuity $\dot{s}$.

\subsection{Variational formulations}

An interesting feature of the equations above is that they can be derived from a variational principle. 
Indeed, the (so-called action) functional $S=\int_{t_1}^{t_2}\int_{x_1}^{x_2}\mathscr{L}\,\ud\/x\,\ud\/t$ 
with the Lagrangian density 
\begin{equation}\label{defL}
\mathscr{L}\ \eqdef\ \rho\,\phi_t\ +\ \half\,\rho\,\phi_x^{\,2}\ +\ \mathscr{V}(\rho),
\end{equation}
where $\mathscr{V}$ is the density of potential energy defined by 
\begin{equation}\label{defV}
\mathscr{V}(\rho)\ \eqdef\ \int\varpi(\rho)\,\ud\/\rho,
\end{equation}
provided that an equation of state $\varpi(\rho)$ (such as (\ref{varpiisen}) given in Appendix \ref{appisent}), 
is substituted into the right-hand side of (\ref{defV}). Since $\varpi(\bar{\rho})=0$ with (\ref{varpiisen}) 
($\bar{\rho}$ a constant state of reference), $\mathscr{V}$ is such that $\mathscr{V}'(\bar{\rho})=0$. 
Note that $\mathscr{V}$ can also be kept explicitly into the Lagrangian if the equation of stated is added 
via a Lagrange multiplier $\lambda$, i.e., considering the Lagrangian density
\begin{align*}
\mathscr{L}'\ \eqdef\ \mathscr{L}\ +\,\left\{\mathscr{V}(\rho)-\int\varpi(\rho)\,\ud\/\rho\right\}\lambda.
\end{align*}  
This is of no interest here, however, so we do not consider this generalisation, for simplicity.
The Euler--Lagrange equations for the Lagrangian density (\ref{defL}) yield 
\begin{align*}
\delta\phi:\quad & 0\ =\ \rho_t\ +\,\left[\,\rho\,\phi_x\,\right]_x,\\
\delta\rho:\quad & 0\ =\ \phi_t\ +\ \half\,\phi_x^{\,2}\ +\ \mathscr{V}'(\rho),
\end{align*}
so the equations of motion (\ref{eqEula}) and (\ref{eqbern}) are recovered. 

An alternative variational formulation is obtained from the Hamilton principle yielding the 
Lagrangian density 
\begin{align}\label{defL0} 
\mathscr{L}_0\ \eqdef\ \half\,\rho\,u^2\ -\ \mathscr{V}(\rho)\  +\,\left\{\,\rho_t\,+ 
\left[\,\rho\,u\,\right]_x\,\right\}\phi,
\end{align}
that is the kinetic minus potential energies plus a constraint enforcing the mass conservation. 
The Euler--Lagrange equations for (\ref{defL0}) yield 
\begin{align*}
\delta\phi:\quad & 0\ =\ \rho_t\ +\,\left[\,\rho\,u\,\right]_x,\\
\delta u:\quad & 0\ =\ u\ -\ \phi_x, \\
\delta\rho:\quad & 0\ =\ \half\,u^2\ -\ \mathscr{V}'(\rho)\ -\ \phi_t\ -\ u\,\phi_x,
\end{align*}
and, substituting $\varpi$ for $\mathscr{V}'$, the barotropic equations (\ref{eqEula}) and 
(\ref{mompi}) are recovered.

The two variational principles above differ by boundary terms only, i.e.,
\begin{equation*}
\mathscr{L}\ +\ \mathscr{L}_0\ =\,\left[\,\rho\,\phi\,\right]_t\ +\,
\left[\,\rho\,u\,\phi\,\right]_x\ +\ \half\,\rho\left(\,u\,-\,\phi_x\,\right)^2, 
\end{equation*}
so the right-hand side yields only boundary terms since $u=\phi_x$.
As advocated by \citet{ClamondDutykh2013}, the variational Hamilton principle is more useful 
for practical applications; this point is illustrated in the section \ref{secregbar} below.

\section{Regularised barotropic flows}\label{secregbar}

Here, we seek for a regularisation of the barotropic Euler equations. We give 
some heuristic arguments for the derivation of such models.

\subsection{Modified Lagrangian}

Following the regularisation of \citet{ClamondDutykh2018a} for the Saint-Venant shallow 
water equations, we seek  for a regularisation of the barotropic Euler equation modifying 
the Lagrangian density as
\begin{align*}
\overline{\mathscr{L}}_\epsilon\ \eqdef\ \mathscr{L}_0\ +\ \epsilon\,\mathscr{A}(\rho)\left[\,
u_t\,+\,u\,u_x\,\right]_x\ +\ \epsilon\,\mathscr{B}(\rho)\left[\,\mathscr{V}''(\rho)\,\rho_x\,
\right]_x,
\end{align*}
where $\epsilon\geqslant0$ is a real parameter at our disposal and $\mathscr{A}$ and $\mathscr{B}$ 
are functions of $\rho$ to be chosen later with suitable properties. 

Note that we could also seek for modifications separating $u_t$ and $uu_x$ in the additional terms 
--- i.e., replacing $\epsilon\mathscr{A}(\rho)\left[u_t+uu_x\right]_x$ by $\epsilon\mathscr{A}(\rho)
u_{xt}+\epsilon\mathscr{C}(\rho)\left[uu_x\right]_x$ --- but that would break the Galilean invariance. 
So, $\mathscr{C}=\mathscr{A}$ is the only physically admissible possibilities.

Exploiting the relations
\begin{align*}
\mathscr{A}(\rho)\left[\,u_t\,+\,u\,u_x\,\right]_x\ &=\,\left[\,\mathscr{A}(\rho)\,u_x\,
\right]_t\ +\ \left[\,\mathscr{A}(\rho)\,u\,u_x\,\right]_x\ +\ \mathscr{A}'(\rho)\,\rho\,
u_x^{\,2},\\
\mathscr{B}(\rho)\left[\,\mathscr{V}''(\rho)\,\rho_x\,\right]_x\ &=\,\left[\,\mathscr{B}(\rho)\,
\mathscr{V}''(\rho)\,\rho_x\,\right]_x\ -\ \mathscr{B}'(\rho)\,\mathscr{V}''(\rho)\,\rho_x^{\,2},
\end{align*}
we derive the equivalent simplified Lagrangian density
\begin{align}\label{defLeps} 
\mathscr{L}_\epsilon\ \eqdef\,\half\,\rho\,u^2\ +\ \epsilon\,\mathscr{A}'\,\rho\,u_x^{\,2}\ 
-\ \mathscr{V}\ -\ \epsilon\,\mathscr{B}'\/\mathscr{V}''\,\rho_x^{\,2}\ +\,\left\{\,\rho_t\,
+\left[\,\rho\,u\,\right]_x\,\right\}\phi.
\end{align}
The functionals given by ${\mathscr{L}}_\epsilon$ and $\overline{\mathscr{L}}_\epsilon$ 
differing only by boundary terms (i.e., $\mathscr{L}_\epsilon-\overline{\mathscr{L}}_\epsilon
=[\cdots]_t+[\cdots]_x$), they yield the same equations of motion.

From (\ref{defLeps}), the regularised kinetic and potential energy densities, respectively 
$\mathscr{K}_\epsilon$ and $\mathscr{V}_\epsilon$, are
\begin{equation*}
\mathscr{K}_\epsilon\ \eqdef\ \half\,\rho\,u^2\ +\ \epsilon\,\mathscr{A}'\,\rho\,u_x^{\,2}, 
\qquad
\mathscr{V}_\epsilon\ \eqdef\ \mathscr{V}\ +\ \epsilon\,\mathscr{B}'\,\mathscr{V}''\,
\rho_x^{\,2}.
\end{equation*}
The total energy is then 
\begin{equation*}
\mathscr{H}_\epsilon\ \eqdef\ \half\,\rho\,u^2\ +\ \epsilon\,\mathscr{A}'\,\rho\,u_x^{\,2}\  
+\ \mathscr{V}\ +\ \epsilon\,\mathscr{B}'\,\mathscr{V}''\,\rho_x^{\,2}.
\end{equation*}
Note that these energies are positive for all $\epsilon\geqslant0$ if $\mathscr{A},\mathscr{B}$ and  
$\mathscr{V}'$ are increasing functions of $\rho$.

\subsection{Linearised equations}\label{:secdispersion}

Here, we consider small perturbations around the rest state $\rho=\bar{\rho}$, $u=0$ and 
$\phi=0$, $\bar{\rho}$ being a positive constant. Introducing $\rho=\bar{\rho}+\tilde{\rho}$, 
$u\approx\tilde{u}$, $\phi\approx\tilde{\phi}$, and $f_0 \eqdef f(\bar{\rho})$ for any function $f$, the tilde quantities being assumed 
small, an approximation of $\mathscr{L}_\epsilon$ up to the second-order is
\begin{align*}
\widetilde{\mathscr{L}}_\epsilon\ =\ \half\,\bar{\rho}\,\tilde{u}^2\ +\ \epsilon\,\mathscr{A}_0'
\,\bar{\rho}\,\tilde{u}_x^{\,2}\ -\ \mathscr{V}_0\ -\ \half\mathscr{V}_0''\,\tilde{\rho}^2\ -\ 
\epsilon\,\mathscr{B}_0'\/\mathscr{V}_0''\,\tilde{\rho}_x^{\,2}\ +\,\left(\tilde{\rho}_t
+\bar{\rho}\tilde{u}_x\right)\tilde{\phi}.
\end{align*}
The Euler--Lagrange (linear) equations for this approximate Lagrangian are
\begin{align}
\delta\tilde{\phi}:\quad & 0\ =\ \tilde{\rho}_t\ +\ \bar{\rho}\,\tilde{u}_x,\label{dL2dphi}\\
\delta\tilde{u}:\quad & 0\ =\ \tilde{u}\ -\ 2\,\epsilon\,\mathscr{A}_0'\,\tilde{u}_{xx}\ 
-\ \tilde{\phi}_x, \label{dL2du}\\
\delta\tilde{\rho}:\quad & 0\ =\ \mathscr{V}_0''\,\tilde{\rho}\ -\ 
2\,\epsilon\,\mathscr{B}_0'\/\mathscr{V}_0''\,\tilde{\rho}_{xx}\ +\ \tilde{\phi}_t.\label{dL2drho}
\end{align}
Looking for traveling waves of the form $\tilde{\rho}=R\cos(kx-\omega t)$, 
$\tilde{u}=U\cos(kx-\omega t)$ and $\tilde{\phi}=\Phi\sin(kx-\omega t)$,  
the equations (\ref{dL2dphi})--(\ref{dL2drho}) yield $\Phi=(1+2\epsilon k^2\mathscr{A}_0')U/k$, 
$U=\omega R/k\bar{\rho}$ and the dispersion relation
\begin{equation*}
\frac{\omega^2}{k^2}\ =\ \bar{\rho}\,\mathscr{V}_0''\,\frac{1\,+\,2\,\epsilon\,k^2
\,\mathscr{B}_0'}{1\,+\,2\,\epsilon\,k^2\,\mathscr{A}_0'}.
\end{equation*}
If $\epsilon=0$ the wave is dispersionless, i.e., the phase velocity $c\eqdef\omega/k$ is 
independent of the wave number $k$. If $\epsilon>0$, the wave is dispersionless if 
$\mathscr{B}_0'=\mathscr{A}_0'$. This condition should be satisfied for all $\bar{\rho}$ and 
for all possible (barotropic) equation of state. Thus, we should take 
\begin{equation}\label{defAB}
\mathscr{B}(\rho)\ =\ \mathscr{A}(\rho).
\end{equation}
Hereafter, we consider only the special case (\ref{defAB}) because we are only interested 
by non-dispersive regularisations of the barotropic Euler equations.

\subsection{Equations of motion}

With (\ref{defAB}), the Euler--Lagrange equations for the Lagrangian density 
(\ref{defLeps}) yield
\begin{align}
\delta\phi:\quad & 0\ =\ \rho_t\ +\,\left[\,\rho\,u\,\right]_x,\label{dLepsdphi}\\
\delta u:\quad & 0\ =\ \rho\,u\ -\ 2\,\epsilon\left[\,\mathscr{A}'\,\rho\,u_x
\,\right]_x\  -\ \rho\,\phi_x, \label{dLepsdu}\\
\delta\rho:\quad & 0\ =\ \half\,u^2\ +\ \epsilon\,(\mathscr{A}'+\mathscr{A}''\rho)\,u_x^{\,2}\ 
-\ \mathscr{V}'\ +\ \epsilon\left(\mathscr{A}''\mathscr{V}''+\mathscr{A}'\/\mathscr{V}'''\right)
\rho_x^{\,2}\nonumber\\
&\qquad +\ 2\,\epsilon\,\mathscr{A}'\,\mathscr{V}''\,\rho_{xx}\ -\ 
\phi_t\ -\ u\,\phi_x, \label{dLepsdrho}
\end{align}
thence
\begin{align*}
\phi_x\ &=\ u\ -\ 2\,\epsilon\,\rho^{-1}\left[\,\mathscr{A}'\,\rho\,u_x
\,\right]_x,\\
\phi_t\ &=\ -\half\,u^2\ +\ \epsilon\,(\mathscr{A}'+\mathscr{A}''\rho)\,u_x^{\,2}\ 
-\ \mathscr{V}'\ +\ \epsilon\left(\mathscr{A}''\mathscr{V}''+\mathscr{A}'\/\mathscr{V}'''\right)
\rho_x^{\,2}\nonumber\\
&\qquad +\ 2\,\epsilon\,\mathscr{A}'\,\mathscr{V}''\,\rho_{xx}\ +\ 2\,\epsilon\,u\,\rho^{-1}
\left[\,\mathscr{A}'\,\rho\,u_x\,\right]_x.  
\end{align*}
Eliminating $\phim$ between these last two relations one obtains
\begin{gather}
0\ =\ \partial_t\!\left\{\,u\, -\, 2\,\epsilon\,\rho^{-1}\left[\,\mathscr{A}'\,\rho\,u_x
\,\right]_x\,\right\}\,+\ 
\partial_x\!\left\{\,\half\,u^2\, +\, \mathscr{V}'\,-\, 2\,\epsilon\,\mathscr{A}'\,\mathscr{V}''
\,\rho_{xx}\right.\nonumber\\
\left.\, -\, \epsilon\,(\mathscr{A}'+\mathscr{A}''
\rho)\,u_x^{\,2}\,-\, \epsilon\left(\mathscr{A}''\mathscr{V}''+\mathscr{A}'\/\mathscr{V}'''\right)
\rho_x^{\,2}\,-\,2\,\epsilon\,u\,\rho^{-1}
\left[\,\mathscr{A}'\,\rho\,u_x\,\right]_x\,\right\}. \label{momtan}
\end{gather}
The equations (\ref{dLepsdphi}) and (\ref{momtan}) form the regularised Euler equations 
for barotropic motions studied in the present paper.

\subsection{Secondary equations}

From the regularised barotropic Euler (rbE) equations (\ref{dLepsdphi}) and (\ref{momtan}), 
several secondary equations can be derived; in particular:
\begin{align}
u_t\ +\ u\,u_x\ +\ \varpi_x\ +\ \epsilon\,\rho^{-1}\,\mathscr{R}_x\, &=\ 0, \label{eqmomnc}\\
\left[\,\rho\,u\,\right]_t\ +\,\left[\,\rho\,u^2\,+\,\rho\,\mathscr{V}'\,-\,\mathscr{V}\,
+\,\epsilon\,\mathscr{R}\,\right]_x\, &=\ 0, \label{eqmomflu}\\
m_t\, +\left[\,u\,m\,+\,\rho\/\mathscr{V}'\,-\,\mathscr{V}\,-\,\epsilon\left(\rho^2\/
\mathscr{A}'\right)'\/u_x^{\,2}\,-\,2\,\epsilon\,\rho\/\mathscr{A}'\,\varpi_{xx}\,+\,
\epsilon\left(\/\rho\/\mathscr{V}''\//\/\mathscr{A}'\right)'\/\mathscr{A}_x^{\,2}\,\right]_x\, &=\ 0, 
\label{eqmomm}\\
\left[\,\half\,\rho\,u^2\,+\,\epsilon\,\rho\,\mathscr{A}'\,u_x^{\,2}\,+\,\mathscr{V}\,+\,
\epsilon\,\mathscr{A}'\,\mathscr{V}''\,\rho_x^{\,2}\,\right]_t\ +\qquad &\nonumber\\
\left[\/\left(\,\half\,\rho\,u^2\,+\,\rho\,\mathscr{V}'\,+\,\epsilon\,\rho\,\mathscr{A}'\,u_x^{\,2}\,
+\,\epsilon\,\mathscr{A}'\,\mathscr{V}''\,\rho_x^{\,2}\,+\,\epsilon\,\mathscr{R}\,\right)u\,
+\,2\,\epsilon\,\rho\,\mathscr{A}'\,\mathscr{V}''\,\rho_x\,u_x\,\right]_x\, &=\ 0, \label{eqene}
\end{align}
where
\begin{align}
\mathscr{R}\ &\eqdef\,\left(\/\rho^2\/\mathscr{A}'\/\right)'u_x^{\,2}\ -\ 2\,\rho\,\mathscr{A}'
\left[\,u_t\,+\,u\,u_x\,+\,\varpi_x\,\right]_x\ +\ \,\left(\/\rho\/\mathscr{V}''\/
/\/\mathscr{A}'\right)'\,\mathscr{A}_x^{\,2}, \label{defR}\\
m\ &\eqdef\ \rho\,u\ -\ 2\,\epsilon\left[\,\rho\,\mathscr{A}'\,u_x\,\right]_x. \nonumber
\end{align}

Introducing the linear Sturm--Liouville operator $\mathcal{L}_\rho\eqdef\/\rho\/-\/
2\/\epsilon\/\partial_x\/\rho\/\mathscr{A}'\/\partial_x$, the equation (\ref{eqmomnc}) multiplied by 
$\rho$ becomes
\begin{equation*}
\mathcal{L}_\rho\!\left\{\,u_t\,+\,u\,u_x\,+\,\varpi_x\,\right\}\, +\ \epsilon\left[\,\left(\/\rho^2\/
\mathscr{A}'\/\right)'u_x^{\,2}\, +\left(\/\rho\/\mathscr{V}''\//\/\mathscr{A}'\/\right)'
\,\mathscr{A}_x^{\,2}\,\right]_x\ =\ 0, 
\end{equation*}
or, inverting the operator,
\begin{equation}\label{momG}
u_t\ +\ u\,u_x\ +\ \varpi_x\ =\ -\,\epsilon\,\mathcal{G}_\rho\partial_x\left\{\,\left(\/\rho^2\/
\mathscr{A}'\/\right)'u_x^{\,2}\, +\left(\/\rho\/\mathscr{V}''\//\/\mathscr{A}'\right)'
\,\mathscr{A}_x^{\,2}\,\right\}, 
\end{equation}
where $\mathcal{G}_\rho=\mathcal{L}_\rho^{-1}$. The operator $\mathcal{G}_\rho\/\partial_x$ 
acting on high frequencies like a first-order anti-derivative, it has a smoothing effect. 
However, this equation is in a non-conservative form. A conservative variant is obtained 
multiplying (\ref{momG}) by $\rho$ and exploiting the mass conservation, hence
\begin{equation}
\left[\,\rho\,u\,\right]_t\ +\,\left[\,\rho\,u^2\,+\,\rho\,\mathscr{V}'\,-\,\mathscr{V}\,
+\,\epsilon\,\mathcal{J}_\rho\!\left\{\,\left(\/\rho^2\/\mathscr{A}'\/\right)'u_x^{\,2}\, 
+\left(\/\rho\/\mathscr{V}''\//\/\mathscr{A}'\right)'\,\mathscr{A}_x^{\,2}\,\right\}\,
\right]_x\ =\ 0, \label{eqmomfluD}
\end{equation}
with the operator
\begin{align*}
\mathcal{J}_\rho\ &\eqdef\ \partial_x^{-1}\,\rho\,\mathcal{G}_\rho\,\partial_x\ =\ 
\partial_x^{-1}\,\rho\left[\,1\,-\,2\,\epsilon\,\rho^{-1}\,\partial_x\,\rho\,
\mathscr{A}'\,\partial_x\,\right]^{-1}\rho^{-1}\,\partial_x \nonumber\\
&=\,\left[\,1\,-\,2\,\epsilon\,\rho\,\mathscr{A}'\,\partial_x\,\rho^{-1}\,\partial_x\,\right]^{-1} 
\ =\ 1\ +\ 2\, \varepsilon\, \rho\, \mathscr{A}'\, \partial_x\, \mathcal{G}_\rho\, \partial_x.
\end{align*}
Comparing (\ref{eqmomfluD}) with (\ref{eqmomflu}), one obtains at once an alternative 
expression for the regularising term 
\begin{align}
\mathscr{R}\ =\ \mathcal{J}_\rho\!\left\{\,\left(\/\rho^2\/\mathscr{A}'\/\right)'u_x^{\,2}\, 
+\left(\/\rho\/\mathscr{V}''\//\/\mathscr{A}'\right)'\,\mathscr{A}_x^{\,2}\,\right\}. \label{defRbis}
\end{align}
While the definition (\ref{defR}) of $\mathscr{R}$ involves second-order spacial derivatives, the 
alternative form (\ref{defRbis}) shows actually that $\mathscr{R}$ behaves at high frequencies somehow like 
zeroth-order derivatives.  Moreover, since the relations (\ref{defR}) and (\ref{defRbis}) are identical, 
we obtain yet another form of the momentum equation
\begin{equation}\label{eqmomter}
2\,\rho\,\mathscr{A}'\left[\,u_t\,+\,u\,u_x\,+\,\varpi_x\,\right]_x\ +\, \left(\mathcal{J}_\rho
-\mathcal{I}\right)\!\left\{\,\left(\/\rho^2\/\mathscr{A}'\/\right)'u_x^{\,2}\, 
+\left(\/\rho\/\mathscr{V}''\//\/\mathscr{A}'\right)'\,\mathscr{A}_x^{\,2}\,\right\}\,=\ 0,
\end{equation}
where $\mathcal{I}$ is the identity operator. Note that applying the operator $\partial_x^{-1}
(2\/\rho\/\mathscr{A}')^{-1}$, the equation (\ref{eqmomter}) can be rewritten
\begin{equation*}
u_t\ +\ u\,u_x\ +\ \varpi_x\ +\ \partial_x^{-1}\,(2\,\rho\,\mathscr{A}')^{-1}
\left(\mathcal{I}-\mathcal{J}_\rho^{-1}\right)\!\{\mathscr{R}\}\ =\ 0,
\end{equation*}
and with
\begin{equation*}
\mathcal{I}\ -\ \mathcal{J}_\rho^{-1}\ =\ \mathcal{I}\ -\ \partial_x^{-1}\,\rho\left[\,1\,
-\,2\,\epsilon\,\rho^{-1}\,\partial_x\,\rho\,\mathscr{A}'\,\partial_x\,\right]\rho^{-1}\,
\partial_x\ =\ 2\,\epsilon\,\rho\,\mathscr{A}'\,\partial_x\,\rho^{-1}\,\partial_x,
\end{equation*}
one gets the equation (\ref{eqmomnc}), as it should be.

\subsection{Rankine--Hugoniot conditions}

Here, we assume that $\rho_x$ and $u_x$ are both continuous if $\epsilon>0$ and that discontinuities 
(if any) occur only in $\rho_{xx}$ and $u_{xx}$. Differentiating twice with respect of $x$ the mass 
conservation (\ref{dLepsdphi}), the jump condition of the resulting equation is 
\begin{equation*}
\left(u\,-\,\dot{s}\right)\llbracket\,\rho_{xx}\,\rrbracket\,+\ \rho\,\llbracket\,u_{xx}\,
\rrbracket\, =\ 0,
\end{equation*}
while the jump condition for (\ref{momtan}) is (provided that $\epsilon$ and $\mathscr{A}'$ 
are not zero)
\begin{equation*}
\left(u\,-\,\dot{s}\right)\llbracket\,u_{xx}\,\rrbracket\,+\ 
\mathscr{V}''\,\llbracket\,\rho_{xx}\,\rrbracket\, =\ 0.
\end{equation*} 
Thus, the speed of the singularity is identical to the original one \eqref{velocityeuler}, whatever the function 
$\mathscr{A}'\neq0$ is. Therefore, a suitable choice for the function $\mathscr{A}$ cannot be 
determined by this consideration.

\subsection{Hamiltonian formulation}
Following the Hamiltonian formulation of the Green--Naghdi equations proposed in \cite{Li}, we introduce the momentum $m\eqdef\rho u-2\epsilon\left[\/\rho\/\mathscr{A}'\/u_x\/\right]_x$ 
and the Hamiltonian functional density 
\begin{equation*} 
\mathscr{H}_\epsilon(\rho,m)\ \eqdef\ \half\,m\,\mathcal{G}_\rho\{m\}\ +\ \mathscr{V}\ +\ 
\epsilon\,\mathscr{A}'\,\mathscr{V}''\,\rho_x^{\,2},
\end{equation*}
we have
\begin{align*}
\mathcal{E}_m\{\mathscr{H}_\epsilon\}\ &=\ \mathcal{G}_\rho\{m\}\ =\ u, \\
\mathcal{E}_\rho\{\mathscr{H}_\epsilon\}\ &=\ \mathscr{V}'\ -\ \epsilon\left(\mathscr{A}'\mathscr{V}''
\right)'\,\rho_x^{\,2}\ -\ 2\,\epsilon\,\mathscr{A}'\mathscr{V}''\,\rho_{xx}\ -\ \half\,u^2\ -\ \epsilon\left(\rho\mathscr{A}'\right)'\,u_x^{\,2}, 
\end{align*} 
where $\mathcal{E}_m$ and $\mathcal{E}_\rho$ are the Euler--Lagrange operators with respect 
of $m$ and $\rho$.  
The rbE equations have then the Hamiltonian structure
\begin{gather*}
\partial_t\binom{\rho}{m}\ =\ -\,\mathbb{J}\boldsymbol{\cdot}
\binom{\mathcal{E}_\rho\{\mathscr{H}_\epsilon\}}{\mathcal{E}_m\{\mathscr{H}_\epsilon\}},\qquad 
\mathbb{J}\ \eqdef\ 
\left[
\begin{array}{cc}
0  &  \partial_x\,\rho \\
\rho\,\partial_x  &  m\,\partial_x\,+\,\partial_x\,m    
\end{array}
\right],
\end{gather*}
yielding the equations (\ref{dLepsdphi}) and (\ref{eqmomm}). It should be noted that 
$\mathbb{J}$ being skew-symmetric and satisfying the Jacobi identity \citep{Nutku1987}, it 
is a proper Hamiltonian (Lie--Poisson) operator.

\subsection{Steady motions}

We seek here for solutions independent of the time $t$, i.e., we look for travelling waves 
of permanent form observed in the frame of reference moving with the wave (note that the 
rbE equations are Galilean invariant). For such flows, the mass conversation yields
\begin{equation}\label{defI}
u\ =\ I\,/\,\rho,
\end{equation}
where $I$ is an integration constant (the mean impulse). From the relations (\ref{eqmomflu}) 
and (\ref{eqene}), the mean (constant) momentum and energy fluxes are respectively
\begin{align}
S\ &=\ \rho\,u^2\,+\,\rho\,\mathscr{V}'\,-\,\mathscr{V}\,
+\,\epsilon\,\mathscr{R}, \label{defS} \\ 
F\ &=\ \left(\,\half\,\rho\,u^2\,+\,\rho\,\mathscr{V}'\,+\,\epsilon\,\rho\,\mathscr{A}'\,u_x^{\,2}\,
+\,\epsilon\,\mathscr{A}'\,\mathscr{V}''\,\rho_x^{\,2}\,+\,\epsilon\,\mathscr{R}\,\right)u\,
+\,2\,\epsilon\,\rho\,\mathscr{A}'\,\mathscr{V}''\,\rho_x\,u_x, \label{defF}
\end{align}
thence --- eliminating $\mathscr{R}$ and using (\ref{defI}) --- the ordinary differential equation
\begin{equation}\label{ODE}
\frac{2\,\epsilon\,\mathscr{A}'}{\rho^2}\/
\left(\frac{\ud\,\rho}{\ud\/x}\right)^{\!2}\ =\ \frac{I^2\,-\,2\,S\,\rho\,+\,2\,(F/I)\,\rho^2
\,-\,2\,\rho\,\mathscr{V}}{I^2\,-\,\rho^3\,\mathscr{V}''}.
\end{equation}

Considering equilibrium states in the far field --- i.e., $\rho\to\rho_\pm$ and $u\to u_\pm$ 
as $x\to\pm\infty$, $\rho_\pm$ and $u_\pm$ being constants --- we have $\mathscr{R}\to0$ and 
the fluxes in the far field are
\begin{align} \label{farflux}
I_\pm\ \eqdef\ \rho_\pm\,u_\pm, \quad 
S_\pm\ \eqdef\ \rho_\pm\,u_\pm^{\,2}\ +\ \rho_\pm\,\mathscr{V}_\pm'\ -\ \mathscr{V}_\pm, \quad
F_\pm\ \eqdef\ \half\,\rho_\pm\,u_\pm^{\,3}\ +\ \rho_\pm\,\mathscr{V}_\pm'\,u_\pm.
\end{align}
For regular solutions, the fluxes of mass, momentum and energy are constants, so $I_+=I_-=I$, 
$S_+=S_-=S$ and $F_+=F_-=F$. For weak solutions, however, we assume that only the mass and 
momentum are conserved (i.e., $I_+=I_-=I$ and $S_+=S_-=S$), some energy being lost at the 
singularity (shock) so $F_+\neq F_-$. 

It should be noted that $\mathscr{A}$ does not appear in the relations (\ref{farflux}). 
The role of $\mathscr{A}$ is to control the singularity at the shock. So, {\em a priori}, 
a local analysis of the singularity is necessary to obtain further informations on $\mathscr{A}$.

\subsection{Local analysis of steady solution}\label{locana}

Let assume that we have a (weak) steady solution with far field conditions (\ref{farflux}) 
and with, possibly, only one singularity at $x=0$ where the density is assumed on the form
\begin{equation}\label{rholoc}
\rho\ =\ \bar{\rho}\ +\ \varrho_\pm\,|x|^\alpha\ +\ \mathrm{o}\!\left(|x|^\alpha\right),
\end{equation}
where $\alpha>0$ is a constant to be found. The plus and minus subscripts in $\varrho$ denote 
$x>0$ and $x<0$, respectively. 
With $f_0 \eqdef f(\bar{\rho})$ for any function $f$, the constant mass flux (\ref{defI}) yields
\begin{equation}\label{uloc}
u\ =\ \frac{I}{\bar{\rho}}\left(\,1\, -\ \frac{\varrho_\pm}{\bar{\rho}}\,|x|^\alpha\,\right)\, 
+\ \mathrm{o}\!\left(|x|^\alpha\right),
\end{equation}
thence
\begin{equation}\label{locR} 
\mathscr{R}\ =\ 2\,\alpha\,(\alpha-1)\,\varrho_\pm\,\bar{\rho}^{\,-2}\,\mathscr{A}_0'\left(\,I^2\,-\,
\bar{\rho}^{\,3}\,\mathscr{V}_0''\,\right)|x|^{\alpha-2}\ +\ \mathrm{o}\!\left(|x|^{\alpha-2}\right),
\end{equation}
and the ODE (\ref{ODE}) yields
\begin{equation}\label{ODE0}
\frac{2\/\epsilon\/\mathscr{A}_0'\/\alpha^2\/\varrho_\pm^{\,2}}{\bar{\rho}^{\,2}}\,|x|^{2\alpha-2}\, 
+\, \mathrm{o}\!\left(|x|^{2\alpha-2}\right)=\, \frac{I^2\,-\,2\/S\/\bar{\rho}\,+\,2\/(F_\pm/I)\/
\bar{\rho}^{\,2}\,-\,2\/\bar{\rho}\/\mathscr{V}_0\,+\, \mathrm{O}\!\left(|x|^\alpha\right)}
{I^2\,-\,\bar{\rho}^{\,3}\/\mathscr{V}_0''\, -\, \left(3\/\bar{\rho}^2\/\mathscr{V}''_0\/+\/
\bar{\rho}^3\/\mathscr{V}_0''' \right)\varrho_\pm\, |x|^\alpha  +\, \mathrm{o}\!\left(|x|^\alpha\right)}.
\end{equation}
From (\ref{locR}) there are three (necessary) possibilities to obtain admissible solutions: $\alpha = 1$ or $\alpha>1$ or $I^2=\bar{\rho}^{\,3}\/\mathscr{V}_0''$. 

If $I^2\neq\bar{\rho}^{\,3}\/\mathscr{V}_0''$, the expansions (\ref{rholoc}) and 
(\ref{uloc}) substituted into (\ref{defS}) and (\ref{defF}) show that $S$ and $F$ cannot be 
constant. Therefore, there are no solutions behaving like  (\ref{rholoc}) if $I^2\neq\bar{\rho}^{\,3}
\/\mathscr{V}_0''$.

%
%

If $I^2=\bar{\rho}^{\,3}\/\mathscr{V}_0''$, the equation (\ref{ODE0}) implies that $\alpha=2/3$ if 
$\mathscr{V}_0\neq I^2/2\bar{\rho}\/-\/S\/+\/(F_\pm/I)\/\bar{\rho}$ and $\alpha=1$ if $\mathscr{V}_0=
I^2/2\bar{\rho}\/-\/S\/+\/(F_\pm/I)\/\bar{\rho}$. The latter case does not yield constant $S$ and $F$,  
so it must be rejected. Finally, the only possibility is $\alpha=2/3$ and (\ref{ODE0}) gives
\begin{equation*}
\frac{8\,\epsilon\,\mathscr{A}_0'\,\varrho_\pm^{\,2}}{9\,\bar{\rho}^{\,3}}\ =\ -\/\frac{I^2\,-\,
2\,S\,\bar{\rho}\,+\,2\,(F_\pm/I)\,\bar{\rho}^{\,2}\,-\,2\,\bar{\rho}\,\mathscr{V}_0}{3\,I^2\,+\,\bar{\rho}^{\,4}
\,\mathscr{V}_0'''}.
\end{equation*}

In summary, the local analysis does not gives hints for a suitable choice of $\mathscr{A}$. 
However, as in \cite{pu2018}, we found the interesting feature that stationary weak solutions have universal 
singularities as $|x|^{2/3}$, whatever the potential $\mathscr{V}$ is and for all possible 
regularising functions $\mathscr{A}$. Note that the analysis above does not rule out the 
possibility of different type of singularities such as $|x|^\alpha\/(\log\!|x|)^\beta$.

\section{Local well-posedness of the regularised barotropic Euler system}\label{secapp} 

The aim of this section is to prove the local well-posedness of the regularised barotropic Euler system introduced in Section \ref{secregbar}.

Let $P=P(\rho)$ denotes the pressure, and let $\rho=\tilde{\rho}+\bar{\rho}$, where $\bar{\rho}$ is a positive constant. Let also 
\begin{equation}\label{defV2}
\varpi(\rho)\ \eqdef\ \int_{\bar{\rho}}^\rho \frac{P'(\alpha)}{\alpha}\, \mathrm{d}\alpha, \qquad \mathscr{V} \eqdef\ \int_{\bar{\rho}}^\rho \varpi(\alpha)\, \mathrm{d}\alpha ,
\end{equation}
where the prime denotes the derivative with respect to $\rho$.

Recalling the operator $\mathcal{L}_\rho\eqdef\/\rho\/-\/2\/\epsilon\/\partial_x\/\rho
\/\mathscr{A}'\/\partial_x$ and the system \eqref{dLepsdphi}, \eqref{momG}
\begin{align} 
\rho_t\ +\,\left[\,\rho\,u\,\right]_x\ &=\ 0, \label{mc}\\
u_t\ +\ u\,u_x\ +\ \varpi_x\ &=\ -\,\epsilon\,\mathcal{L}^{-1}_\rho\partial_x \left\{\,\left(\/\rho^2\/\mathscr{A}'\/\right)'u_x^{\,2}\, 
+\left(\/\rho\/\mathscr{V}''\//\/\mathscr{A}'\right)'\,{\mathscr{A}'}^{\/2}\, \rho_x^{\,2}\,\right\},\label{momG2} 
\end{align}
where smooth solutions of \eqref{mc}, \eqref{momG2} satisfy the energy equation \eqref{eqene} with $\mathscr{R}$ is defined in \eqref{defRbis}.
The goal of this section is to prove the following theorem
\begin{thm}\label{existencerEuler}
Let $\tilde{m} \geqslant s \geqslant 2$, $\tilde{m}$ be an integer, $P, \mathscr{A} \in {C}^{\tilde{m}+4}(]0,+\infty[)$ such that $P'(\rho)~>~0$, $ \mathscr{A}'(\rho)>0$ for $\rho>0$. Let also $W_0=(\tilde{\rho}_0,u_0)^\top \in H^s$ satisfying $\inf_{x \in \mathds{R}} \rho_0(x) >~\rho^{*}$, then there exist $T>0$ and a unique solution $W \in {C}([0,T], H^s) \cap {C}^1([0,T], H^{s-1})$ of \eqref{mc}, \eqref{momG2} satisfying the non-emptiness condition 
$ \inf_{x \in \mathds{R}}\, \rho(t,x)\ >\ 0, $ and the conservation of the energy
\begin{equation}\label{energycons}
\frac{\ud}{\ud\/t} \int_\mathds{R} \left(\,\half\,\rho\,u^2\,+\,\epsilon\,\rho\,\mathscr{A}'\,u_x^{\,2}\,+\,\mathscr{V}\,+\,
\epsilon\,\mathscr{A}'\,\mathscr{V}''\,\rho_x^{\,2}\,\right)\, \ud x\ =\ 0.
\end{equation}
Moreover, if the maximal existence time $T_{max}< +\infty$, then 
\begin{equation}\label{bucrE}
\lim_{t \to T_{max}} \|W_x\|_{L^\infty}\ =\ +\infty.
\end{equation}
\end{thm}
\begin{rem}
The solution given in the previous theorem depends continuously on the initial data in the sense: If $W_0, \tilde{W}_0 \in H^s$, such that $\rho_0,\tilde{\rho}_0 \geqslant \rho^*$, then there exists a constant $C(\|\tilde{W}\|_{L^\infty([0,T],H^{s})},\|W\|_{L^\infty([0,T],H^{s})})>0$, such that 
\begin{equation*}
\|W-\tilde{W}\|_{L^\infty([0,T],H^{s-1})}\ \leqslant\ C\, \|W_0-\tilde{W}_0\|_{H^s}.
\end{equation*}
\end{rem}
\begin{rem}
Theorem \ref{existencerEuler} holds also for periodic domains.
\end{rem}

\begin{rem}
From the definition \eqref{defV2} we have $\mathscr{V}''(\rho) = P'(\rho)/\rho$, then, the conditions $P'>0$ and $\mathscr{A}'>0$ are enough to obtain a positive energy in \eqref{energycons}.
\end{rem}

\begin{rem} Note that if $\rho \in [\rho_{\mathrm{inf}}, \rho_{\mathrm{sup}}] \subset ]0,+\infty[$, then $0<\alpha \leqslant P'(\rho)/\rho \leqslant \beta < +\infty$. This implies with the definition \eqref{defV2} that $\alpha (\rho-\bar{\rho})^2 \leqslant \mathscr{V} \leqslant \beta (\rho-\bar{\rho})^2$. Then, the conserved energy \eqref{energycons} is equivalent to the $H^1$ norm of $(\tilde{\rho},u)$.
\end{rem}

\subsection{Preliminary results}
Let $\Lambda$ be defined such that $\widehat{\Lambda f}=(1+\xi^2)^\frac{1}{2} \hat{f}$.
In order to prove Theorem \ref{existencerEuler}, we recall the classical lemmas.
\begin{lem}(\cite{kato1988commutator})
Let $[A,B] \eqdef AB-BA$ be the commutator of the operators $A$ and $B$. If $r\, \geqslant\, 0$, then 
\begin{align}
\|f\, g\|_{H^r}\ &\lesssim\ \|f\|_{L^\infty}\, \|g\|_{H^r}\ +\ \|f\|_{H^r}\, \|g\|_{L^\infty}, \label{Algebra} \\
\left\| \left[ \Lambda^r,\, f \right]\, g  \right\|_{L^2}\ &\lesssim\  \|f_x\|_{L^\infty}\, \|g\|_{H^{r-1}}\ +\ \|f\|_{H^r}\, \|g\|_{L^\infty}. \label{Commutator}
\end{align}
\end{lem}

\begin{lem}(\cite{constantin2002initial})
Let $F \in {C}^{\tilde{m}+2}(\mathds{R})$ with $F(0)=0$ and $0 \leqslant s \leqslant \tilde{m}$, then there exists a continuous function $\tilde{F}$, such that for all $f \in H^s \cap W^{1,\infty}$ we have
\begin{equation}\label{Composition2}
\|F(f)\|_{H^s}\ \leqslant\ \tilde{F} \left( \|f\|_{W^{1,\infty}} \right)\, \|f\|_{H^s}.
\end{equation}
\end{lem}

In the following lemma, we prove the invertibility of the operator $\mathcal{L}_\rho$ \eqref{op:SL} and we obtain some estimates satisfied by $\mathcal{L}_\rho^{-1}$.

\begin{lem}\label{Inverseesitimates}
Let $0< \rho_{\mathrm{inf}} \leqslant \rho \in W^{1,\infty}$ and $\mathscr{A} \in {C}^2(]0,+\infty[)$ satisfying $ \mathscr{A}' > 0$, then the operator $\mathcal{L}_\rho$ is an isomorphism from $H^2$ to $L^2$ and 
\begin{enumerate}
\item  If $0 \leqslant s \leqslant \tilde{m} \in \mathds{N}$ and $\mathscr{A} \in {C}^{\tilde{m}+3}(]0,+\infty[)$, then  
\begin{subequations}\label{estimate1}
\begin{align}
\left\|\mathcal{L}_\rho^{-1}\, \partial_x\, \psi \right\|_{H^{s+1}}\  
&\lesssim\ \left\|\psi \right\|_{H^s}\ +\ \left\|\rho\, -\, \bar{\rho} \right\|_{H^s}\, 
\left\|\mathcal{L}_\rho^{-1}\, \partial_x\, \psi \right\|_{W^{1,\infty}},\\
\left\|\mathcal{L}_\rho^{-1}\, \phi \right\|_{H^{s+1}}\  
&\lesssim\ \left\|\phi \right\|_{H^s}\ +\ \left\|\rho\, -\, \bar{\rho} \right\|_{H^s}\, 
\left\|\mathcal{L}_\rho^{-1}\, \phi \right\|_{W^{1,\infty}}.
\end{align}
\end{subequations}
\item  If $0 \leqslant s \leqslant \tilde{m} \in \mathds{N}$ and $\mathscr{A} \in {C}^{\tilde{m}+3}(]0,+\infty[)$, then 
\begin{subequations}\label{estimate2}
\begin{align}
\left\|\mathcal{L}_\rho^{-1}\, \partial_x\, \psi \right\|_{H^{s+1}}\  
&\lesssim\ \left\|\psi \right\|_{H^s}\, \left(1 +\ \left\|\rho\, -\, \bar{\rho} \right\|_{H^s} \right), \\
\left\|\mathcal{L}_\rho^{-1}\, \phi \right\|_{H^{s+1}}\  
&\lesssim\ \left\|\phi \right\|_{H^s}\, \left(1 +\ \left\|\rho\, -\, \bar{\rho} \right\|_{H^s} \right), 
\end{align}
\end{subequations}
\item If $\phi \in {C}_{\mathrm{lim}} \eqdef \{f\in {C},\ f(\pm \infty) \in \mathds{R} \}$, then $\mathcal{L}_\rho^{-1}\, \phi$ is well defined and 
\begin{equation}\label{estimate3}
\left\|\mathcal{L}_\rho^{-1}\, \phi\right\|_{W^{2,\infty}}\ \lesssim\ \left\|\phi\right\|_{L^\infty}. 
\end{equation}
\item If $\psi \in {C}_{\mathrm{lim}} \cap L^1$, then 
\begin{equation}\label{estimate4}
\left\|\mathcal{L}_\rho^{-1}\, \partial_x \psi\right\|_{W^{1,\infty}}\ \lesssim\ \left\|\psi\right\|_{L^\infty}\ +\ \left\|\psi\right\|_{L^1}. 
\end{equation}
\end{enumerate}
All the constants depend on $s, \varepsilon, \rho_\mathrm{inf}, \|\rho - \bar{\rho}\|_{W^{1,\infty}}$ and not on $\|\rho - \bar{\rho}\|_{H^s}$.
\end{lem}
The previous lemma is proven in \cite{liu2019well} for the special case $\mathscr{A}(\rho)=\rho^3/6$. Here, the same proof is followed.

\proof \textbf{Step 0:} In the first step, we prove, using the Lax--Milgram theorem, that $\mathcal{L}_\rho$ is an isomorphism from $H^2$ to $L^2$, let the bi-linear function $a$ from $H^1 \times H^1$ to $\mathds{R}$ such that 
$$ a(u,v)\ \eqdef\ (\rho\, u, v)\ +\ 2\, \varepsilon\, \left( \rho\, \mathscr{A}'\, u_x, v_x \right). $$ 
Using that $\rho$ is bounded and far from zero, one can easily show that the function $a$ is continuous and coercive, then Lax--Milgram theorem shows that there exists a continuous bijection $J$ between $H^1$ and $H^{-1}$, such that for all $u,v \in H^1$ we have 
$$ a(u,v)\ =\ (J u, v)_{H^{-1} \times H^1}.$$  
If $J u \in L^2$, and integration by parts shows that $2 \varepsilon \left( \rho \mathscr{A}' u_x \right)_x = \rho u -J u \in L^2$ and $J=\mathcal{L}_\rho$, this implies that $u \in H^2$ which finishes the proof that $\mathcal{L}_\rho$ is an isomorphism from $H^2$ to $L^2$.

\noindent 
\textbf{Step 1:} Let $\mathcal{L}_\rho u = \phi + \psi_x$, then
\begin{align}
\|u\|_{H^1}^2\ 
&=\ (u,u)\ +\ (u_x,u_x) \nonumber \\ 
&\lesssim\ (\rho\, u,u)\ +\ 2\, \varepsilon\, (\rho\, \mathscr{A}'\, u_x,u_x) \nonumber \\
&=\ (\mathcal{L}_\rho u,u)\ =\ (\phi,u)\ -\ (\psi,u_x) \nonumber\\
&\lesssim\ \|u\|_{H^1}\, \left( \|\phi\|_{L^2}\ +\ \|\psi\|_{L^2} \right), \nonumber
\end{align}
which implies that 
\begin{equation}\label{s=0}
\|u\|_{H^1}\  \lesssim\ \|\phi\|_{L^2}\ +\ \|\psi\|_{L^2}.
\end{equation}
Using the Young inequality $ab \leqslant {\textstyle \frac{1}{2\/ \alpha}} a^2 + {\textstyle \frac{\alpha}{2}} b^2$ with $\alpha>0$ we obtain
\begin{align}
\|u_x\|_{H^1}^2\ 
&=\ (u_x,u_x)\ +\ (u_{xx},u_{xx}) \nonumber \\ 
&\lesssim\ (\rho\, u_x,u_x)\ +\ 2\, \varepsilon\, (\rho\, \mathscr{A}'\, u_{xx},u_{xx}) \nonumber \\
&=\ -(\rho\, u, u_{xx})\ -\ (\rho_x\, u, u_x)\ +\ 2\, \varepsilon\, \left( \left(\rho\, \mathscr{A}' u_x \right)_x\, -\, (\rho\, \mathscr{A}')_x u_x, u_{xx} \right) \nonumber\\
&=\ -(\mathcal{L}_\rho\, u, u_{xx})\ -\ (\rho_x\, u, u_x)\ -\ 2\, \varepsilon\, \left( (\rho\, \mathscr{A}')_x u_x, u_{xx} \right) \nonumber\\
& \lesssim\ \alpha\, \|u_{xx}\|_{L^2}^2\ +\ {\textstyle \frac{1}{\alpha}}\, \left( \|\mathcal{L}_\rho\, u\|_{L^2}^2\ +\ \|u_{x}\|_{L^2}^2 \right)\ +\ \|u\|_{H^1}^2.\nonumber
\end{align}
Taking $\alpha>0$ small enough we obtain that 
\begin{equation*}
\|u_x\|_{H^1}^2\ \lesssim\ \|\mathcal{L}_\rho\, u\|_{L^2}^2\ +\ \|u\|_{H^1}^2.
\end{equation*}
then 
\begin{equation*}
\|u_x\|_{H^1}\ \lesssim\ \|\mathcal{L}_\rho\, u\|_{L^2}\ +\ \|u\|_{H^1}.
\end{equation*}
Taking $\phi=0$ (respectively $\psi=0$) and using \eqref{s=0}, we obtain 
\begin{equation*}
\left\|\mathcal{L}_\rho^{-1}\, \partial_x\, \psi \right\|_{H^2}\ \lesssim\ \|\psi\|_{H^1}\, \qquad \left\|\mathcal{L}_\rho^{-1}\, \phi \right\|_{H^{2}}\  \lesssim\ \|\phi\|_{L^2}.
\end{equation*}
An interpolation with \eqref{s=0} implies that 
\begin{equation}\label{01}
\left\|\mathcal{L}_\rho^{-1}\, \partial_x\, \psi \right\|_{H^{s+1}}\ \lesssim\ \|\psi\|_{H^{s}}, \qquad  
\left\|\mathcal{L}_\rho^{-1}\, \phi \right\|_{H^{s+1}}\ \lesssim\ \|\phi\|_{H^{s}} 
\qquad \forall s \in [0,1].
\end{equation}
Let now $s>0$, and let $\mathcal{L}_\rho u = \phi + \psi_x$, we then have
$$ \mathcal{L}_\rho\, \Lambda^s u\ =\ [\rho, \Lambda^s]\, u\ +\ \Lambda^s \phi\ +\ \partial_x\, \left\{ -\, 2\, \varepsilon\, [\rho\, \mathscr{A}', \Lambda^s]\, u_x\ +\ \Lambda^s \psi \right\}. $$ 
Defining $\tilde{u}=\Lambda^s u$, $\tilde{\phi}=[\rho, \Lambda^s] u +  \Lambda^s \phi $ and $\tilde{\psi}= -2 \varepsilon [\rho \mathscr{A}', \Lambda^s] u_x +  \Lambda^s \psi $ and using \eqref{s=0}, \eqref{Commutator}, \eqref{Composition2} we obtain
\begin{align}
\|\Lambda^s u\|_{H^1}\ 
&\lesssim\ \left\|[\rho, \Lambda^s]\, u \right\|_{L^2}\ +\ \left\|[\rho\, \mathscr{A}', \Lambda^s]\, u_x \right\|_{L^2}\ +\ \left\|\phi \right\|_{H^s}\ +\ \left\|\psi \right\|_{H^s} \nonumber\\
&\lesssim\ \|\Lambda^{s-1} u\|_{H^1}\ +\ \|\rho\, -\, \bar{\rho}\|_{H^s}\, \|u\|_{W^{1,\infty}}\ +\ \left\|\phi \right\|_{H^s}\ +\ \left\|\psi \right\|_{H^s}. \nonumber
\end{align}
Then, by induction (on $s$) and using \eqref{01} one obtains that \eqref{estimate1} holds for all $s \geqslant 0$.

\noindent 
\textbf{Step 2:} If $s \leqslant 1$, then \eqref{estimate2} follows directly from \eqref{01}.
If $s>1$, using the embedding $H^1 \hookrightarrow L^\infty$, \eqref{estimate1} and \eqref{01} for $s=1$ we obtain \eqref{estimate2}.

\noindent 
\textbf{Step 3:} Let $C_0 \eqdef \{f\in {C},\ f(\pm \infty)=0 \}$, using that $L^2 \cap {C}_0$ is dense in ${C}_0$ one can define $\mathcal{L}_\rho^{-1}$ on ${C}_0$.
If $\phi$ is in ${C}_{\mathrm{lim}}$, we use the change of functions (see Lemma 4.4 in \cite{liu2019well}) 
\begin{equation*}
\phi_0(x)\ \eqdef\ \phi(x)\ -\  {\textstyle \frac{1}{\bar{\rho}}}\, \mathcal{L}_\rho \left( \phi(-\infty)\ +\ \left( \phi(+\infty)\ -\ \phi(-\infty) \right)\, \frac{\mathrm{e}^x}{1+\mathrm{e}^x} \right) \in C_0,
\end{equation*}
the operator $\mathcal{L}_\rho^{-1}$ can be defined as
\begin{equation*}
\mathcal{L}_\rho^{-1}\, \phi\ \eqdef\ \mathcal{L}_\rho^{-1}\, \phi_0\ +\  {\textstyle \frac{1}{\bar{\rho}}}\, \left( \phi(-\infty)\ +\ \left( \phi(+\infty)\ -\ \phi(-\infty) \right)\, \frac{\mathrm{e}^x}{1+\mathrm{e}^x} \right).
\end{equation*}
In order order to prove \eqref{estimate3}, let $\phi=\mathcal{L}_\rho u$, using the variable 
\begin{equation}\label{chv}
z\ \eqdef\  \int\/ \frac{\ud\/x}{2\, \rho(x)\, \mathscr{A}'(\rho(x))\,}, 
\end{equation}
we obtain that 
\begin{equation}\label{Helmholtz-like}
\phi\ =\ \rho\, u\ -\  { \frac{\varepsilon}{2\, \rho\, \mathscr{A}'}}\, u_{zz}.
\end{equation}
The classical maximum principle equations implies that $\|u\|_{L^\infty} \leqslant \|\phi\|_{L^\infty}/\rho_{\mathrm{inf}}$, which implies with \eqref{Helmholtz-like} that $\|u_{zz}\|_{L^\infty} \lesssim \|\phi\|_{L^\infty}$, then the Landau--Kolmogorov inequality (see Lemma 4.3 in \cite{liu2019well} for example) implies that $\|u_{z}\|_{L^\infty} \lesssim \|\phi\|_{L^\infty}$. The last inequality with the change of variables \eqref{chv} imply that $\|u_{x}\|_{L^\infty} \lesssim \|\phi\|_{L^\infty}$, using that $2 \varepsilon \rho \mathscr{A}' u_{xx}= \rho u- 2 \varepsilon (\rho \mathscr{A}')_x u_x-\phi$ we obtain \eqref{estimate3}.

\noindent 
\textbf{Step 4:} Note that 
\begin{equation*}
\mathcal{L}_\rho\, \int_{-\infty}^x \frac{\psi}{\rho\, \mathscr{A}'}\, \ud y\ =\ \rho \, \int_{-\infty}^x \frac{\psi}{\rho\, \mathscr{A}'}\, \ud y\ -\ 2\, \varepsilon\, \psi_x.
\end{equation*}
Applying $\mathcal{L}_\rho^{-1}$ and $\partial_x$ one obtains
\begin{align*}
2\, \varepsilon\, \mathcal{L}_\rho^{-1}\, \partial_x\, \psi\ &=\ \mathcal{L}_\rho^{-1}\, \left(\rho\, \int_{-\infty}^x \frac{\psi}{\rho\, \mathscr{A}'}\, \ud y \right)\ -\ \int_{-\infty}^x \frac{\psi}{\rho\, \mathscr{A}'}\, \ud y,\\
2\, \varepsilon\, \partial_x\, \mathcal{L}_\rho^{-1}\, \partial_x\, \psi\ &=\ \partial_x\, \mathcal{L}_\rho^{-1}\, \left(\rho\, \int_{-\infty}^x \frac{\psi}{\rho\, \mathscr{A}'}\, \ud y \right)\ -\ \frac{\psi}{\rho\, \mathscr{A}'}.
\end{align*}
The last two inequalities with \eqref{estimate3} imply \eqref{estimate4}.
\qed

\subsection{Iteration scheme and energy estimate}

Defining $W \eqdef (\tilde{\rho},u)^\top$ and 
\begin{align*}
B(W)\  \eqdef\ \begin{pmatrix} u & \rho \\ \varpi' &  u \end{pmatrix}, \quad
F(W)\ \eqdef\  \begin{pmatrix} 0 \\  -\,\epsilon\,\mathcal{L}^{-1}_\rho\partial_x \left\{\,\left(\/\rho^2\/\mathscr{A}'\/\right)'u_x^{\,2}\, 
+\left(\/\rho\/\mathscr{V}''\//\/\mathscr{A}'\right)'\,{\mathscr{A}'}^{\/2}\, \rho_x^{\,2}\,\right\} \end{pmatrix},
\end{align*}
the system \eqref{mc}, \eqref{momG2} becomes
\begin{equation}\label{generalsystem}
W_t\ +\ B(W)\, W_x\ =\ F(W), \qquad W(0,x)\ =\ W_0(x).
\end{equation}

The proof of the local existence of \eqref{generalsystem} is based on solving the linear hyperbolic system  
\begin{equation}\label{iteration}
\partial_t\, W^{n+1}\ +\ B(W^n)\, \partial_x\, W^{n+1}\ =\ F(W^n), \qquad W^n(0,x)\ =\ (\tilde{\rho}_0(x),u_0(x))^\top,
\end{equation}
for all $n\geqslant0$, where $W^0(t,x)=(\tilde{\rho}_0(x),u_0(x))^\top$. Then, uniform (on $n$) estimates of an 
energy that is equivalent to the $H^s$ norm will be given. Taking the limit $n\to\infty$, we obtain a 
solution of  \eqref{generalsystem}.
Since $\varpi'>0$, the system 
\eqref{iteration} is hyperbolic; which is an important point to solve each iteration in \eqref{iteration}.

For the sake of simplicity, let be $\underline{W} = W^n$ (known on every step of the iteration) and 
let $W=W^{n+1}$ be the solution of the linear system
\begin{equation}\label{linearsystem}
W_{t}\ +\ B(\underline{W})\, W_{x}\ =\ F(\underline{W}), \qquad W(0,x)\ =\ (\tilde{\rho}_0(x),u_0(x))^\top.
\end{equation}
Note that a symmetriser of $\underline{B}=B(\underline{W})$ is 
\begin{equation*}
\underline{A}\ =\ A(\underline{W})\ \eqdef\ \begin{pmatrix} \varpi' & 0 \\ 0 & \rho \end{pmatrix}.
\end{equation*} 
Let the energy of \eqref{linearsystem} be defined as 
\begin{equation*}
E(W)\ \eqdef\ \left(\Lambda^s\, W,\ \underline{A}\, \Lambda^s\, W \right),
\end{equation*}
where $(\cdot,\, \cdot)$ is the scalar product in $L^2$. Since the matrix $\underline{A}\,\underline{B}$ 
is symmetric, a helpful feature for the energy estimates below, it justifies  the use of $\underline{A}$ 
in the definition of the energy $E(W)$.
Note that if $\rho$ is bounded and far from zero, then $E(W)$ is equivalent to $\|W\|_{H^s}^2$.
In order to prove Theorem \ref{existencerEuler}, the following energy estimate is needed.
\begin{thm}\label{energyestimate}
Let $W= (\tilde{\rho},u)^\top$, $\rho=\tilde{\rho}+\tilde{\rho}$, $s \geqslant 2$ and $\rho^{*},R>0$ then there exist $K,T>0$ such that: if the initial data $(\tilde{\rho}_0,u_0) \in H^s$ satisfy
\begin{equation}\label{initial}
\inf_{x \in \mathds{R}} \rho_0(x)\ >\ \rho^{*}, \qquad E(W_0)\ <\ R,
\end{equation}
and $\underline{W} \in {C}([0,T], H^s) \cap {C}^1([0,T], H^{s-1})$, satisfying for all $t \in [0,T]$
\begin{equation*}
\underline{\rho}\ \geqslant\ \rho^{*}, \qquad \|\underline{W}_t\|_{H^{s-1}}\ \leqslant\ K, \qquad E(\underline{W})\ \leqslant\ R,
\end{equation*}
then there exists a unique $W \in C([0,T], H^s) \cap C^1([0,T], H^{s-1})$ a solution of \eqref{linearsystem} such that 
\begin{equation*}
\rho\ \geqslant\ \rho^{*}, \qquad \|W_t\|_{H^{s-1}}\ \leqslant\ K, \qquad E(W)\ \leqslant\ R.
\end{equation*}
\end{thm}

\proof
For the existence of $W$ see Appendix A in \cite{Israwi2011}. Defining $\underline{F}\eqdef F(\underline{W})$ and using \eqref{linearsystem} we obtain 
\begin{align}\label{energyterms}
E(W)_t\ =&\ 2\, \left( \Lambda^s W_t,\, \underline{A}\, \Lambda^s W \right)\ +\ \left( \Lambda^s W,\, \underline{A}_t\, \Lambda^s W \right) \nonumber \\
=&\ -2\, \left( \Lambda^s \underline{B}\, W_x,\, \underline{A}\, \Lambda^s W \right)\ -\ 2\, \left( \Lambda^s \underline{F},\, \underline{A}\, \Lambda^s W \right)\ +\ \left( \Lambda^s W,\, \underline{A}_t\, \Lambda^s W \right) \nonumber\\
=&\ -2\, \left( [\Lambda^s, \underline{B}]\, W_x,\, \underline{A}\, \Lambda^s W \right)\ -\ 2\, \left(\underline{B}\, \Lambda^s  W_x,\, \underline{A}\, \Lambda^s W \right) \nonumber \\
&\ - 2\, \left( \Lambda^s \underline{F},\, \underline{A}\, \Lambda^s W \right)\ +\ \left( \Lambda^s W,\, \underline{A}_t\, \Lambda^s W \right) \nonumber \\
=&:\ I\ +\ II\ +\ III\ +\ IV,
\end{align}
Now, some bounds of the four terms will be given. Note that
\begin{align*}
-\half\, I\, =&\, \left( [\Lambda^s, \underline{u}]\, \tilde{\rho}_x,\, \underline{\varpi'}\,  \Lambda^s \tilde{\rho}  \right)\, 
+\, \left( [\Lambda^s, \underline{\rho}]\, u_x,\, \underline{\varpi'}\, \Lambda^s \tilde{\rho}  \right)\,
+\, \left( [\Lambda^s, \underline{u}]\, u_x,\, \underline{\rho} \Lambda^s u  \right)\,
\\ &
+\, \left( [\Lambda^s, \underline{\varpi'}]\, \tilde{\rho}_x,\, \underline{\rho} \Lambda^s u  \right).
\end{align*}
Using \eqref{Commutator} and \eqref{Composition2} we obtain 
\begin{align*}
\left| \left( [\Lambda^s, \underline{u}]\, \tilde{\rho}_x,\, \underline{\varpi'}\,  \Lambda^s \tilde{\rho}  \right) \right|\ 
&\leqslant\ \left\|[\Lambda^s, \underline{u}]\, \tilde{\rho}_x \right\|_{L^2}\, \left\| \underline{\varpi'}\, \Lambda^s \tilde{\rho} \right\|_{L^2}\\
&\lesssim\  \left( \|\underline{u}\|_{L^\infty}\, \|\tilde{\rho}_x\|_{H^{s-1}}\ +\ \|\underline{u}\|_{H^s}\, \|\tilde{\rho}_x\|_{L^\infty} \right)\, \|\underline{\varpi'} \|_{L^\infty}\,  \left\|\Lambda^s \tilde{\rho} \right\|_{L^2}\\
&\lesssim\ \|W\|_{H^s}^2\ \lesssim\ E(W).
\end{align*}
All the terms of $I$ can be studied by the same way to obtain
\begin{equation}\label{I}
|I|\ \lesssim\  E(W).
\end{equation}
Using that $\underline{A}$ and $\underline{A} \underline{B}$ are symmetric, an integration by parts yield to 
\begin{align}\label{II}
|II|\ =\ \left| \left( \Lambda^s W,\, (\underline{A} \underline{B})_x\, \Lambda^s W  \right)  \right|\
\leqslant\ \| (\underline{A} \underline{B})_x \|_{L^\infty}\, \|W\|_{H^s}^2 \
\lesssim\ \|W\|_{H^s}^2\ \lesssim\ E(W).
\end{align}
Using the Young inequality $2 ab \leqslant a^2+b^2$ one obtains 
\begin{equation*}
|III|\ \leqslant\ \|\underline{A}\|_{L^\infty}\, \left( \|\underline{F}\|_{H^s}^2\ +\ \|W\|_{H^s}^2 \right)
\end{equation*}
From the inequality \eqref{estimate2} we have
\begin{align*}
\|\underline{F}\|_{H^s}\ 
&\lesssim\ \left\|  \left(\/\underline{\rho}^2\/\underline{\mathscr{A}}'\/\right)' \underline{u}_x^{\,2}\, 
+ \left(\/ \underline{\rho}\/\underline{\mathscr{V}}''\//\/\underline{\mathscr{A}}'\right)'\,{\underline{\mathscr{A}}'}^{\/2} \underline{\rho}_x^{\,2} \right\|_{H^{s-1}},
\end{align*}
which implies with \eqref{Algebra} and \eqref{Composition2} that $\|\underline{F}\|_{H^s}$ is bounded, then 
\begin{equation}\label{III}
|III|\ \lesssim\  E(W)\ +\ 1.
\end{equation}
Note that $\left\|\underline{\tilde{\rho}}_t \right\|_{L^\infty} \leqslant \left\|\underline{\tilde{\rho}}_t \right\|_{H^{s-1}} \leqslant K$, then 
\begin{align}\label{IV}
|IV|\ \leqslant\ \|\underline{A}_t\|_{L^\infty}\, \|W\|_{H^s}^2\ \lesssim\ K\, E(W)
\end{align}
The system \eqref{linearsystem} implies
\begin{align}\label{Wt}
\|W_t\|_{H^{s-1}}\ 
&=\ \left\| B(\underline{W})\, W_x\ +\ \underline{F} \right\|_{H^{s-1}}\nonumber\\
&\lesssim\   \|B(\underline{W})\|_{L^\infty}\, \|W_x\|_{H^{s-1}}\ +\ \|B(\underline{W})\|_{H^{s-1}}\, \|W_x\|_{L^\infty}\ +\ \| \underline{F}\|_{H^s} \nonumber \\
&\lesssim\ E(W)\ +\ 1.
\end{align}
All the constants hidden in $"\lesssim"$ do not depend on $K$ and $W$.
Using \eqref{I}, \eqref{II}, \eqref{III} and \eqref{IV} we obtain that 
\begin{equation*}
\partial_t\, E(W)\ \leqslant\ C\, (K+1)\, [ E(W)\ +\ 1],
\end{equation*}
which implies with Gronwall lemma that
\begin{equation}\label{Gronwall}
E(W)\ \leqslant\ [ E(W_0)\ +\ 1]\, \mathrm{e}^{C\, (K+1)\, t}\ -\ 1.
\end{equation}
Since $E(W_0)<R$, choosing first $K>0$ and then $T>0$ such that 
\begin{equation*}
C\, (R\ +\ 1)\ \leqslant\ K, \qquad \qquad [E(W_0)\ +\ 1]\, \mathrm{e}^{C\, (K+1)\, T}\ -\ 1\ \leqslant\ R
\end{equation*}
we obtain with \eqref{Wt} and \eqref{Gronwall} that $\|W_t\|_{H^{s-1}}\ \leqslant\ K$ and $E(W)\ \leqslant\ R$. Since $\|\rho_t\|_{L^\infty} \lesssim \|W\|_{H^{s-1}} \leqslant K$ and $\rho_0>\rho^{*}$ then taking $T$ small enough we have $\rho \geqslant \rho^{*}$. \qed

\subsection{Proof of Theorem \ref{existencerEuler}}

Theorem \ref{energyestimate} shows that if the initial data satisfy \eqref{initial}, then the sequence $(W^n)_{n \in \mathds{N }}$ exists, it is uniformly bounded in $ {C}([0,T], H^s) \cap {C}^1([0,T], H^{s-1})$ and satisfies $\rho^n \geqslant \rho^{*}$. 
Using classical arguments of Sobolev spaces one can prove that there exists $W \in {C}([0,T], H^s)$ such that $W^n$ converges ``up to a sub-sequence'' to $W$ in ${C}([0,T], H^{s'})$ for all $s' \in [0,s[$. Before taking the limit $n \to \infty$ in \eqref{iteration}, we will verify that if $W^n$ converges, then $W^{n+1}$ converges too and towards the same limit. For that purpose, let 
\begin{equation*}
\tilde{E}^n\ \eqdef\ \left( \Lambda^{s-1}\, (W^{n+1}\, -\, W^n),\ A^n\, \Lambda^{s-1}\, (W^{n+1}\, -\, W^n) \right),
\end{equation*}
using estimates as in the proof of Theorem \ref{energyestimate} (see also \cite{majda2012compressible,pu2018} for more details) one can prove that for $T>0$ small enough, we obtain that $\tilde{E}^{n+1} \leqslant  \tilde{E}^n/2$, which implies that $\|W^{n+1}-W^n\|_{L^\infty([0,T],H^{s-1})} \to 0$.
Taking $n \to \infty$ in the weak formulation of  \eqref{iteration}, we obtain that $W$ is a weak solution of \eqref{generalsystem}. Using that $W \in {C}([0,T], H^s)$ and \eqref{generalsystem} we deduce that $W$ is a strong solution and $W \in {C}^1([0,T], H^{s-1})$. 

In order to prove the blow-up criteria \eqref{bucrE}, we suppose that $\|W_x\|_{L^\infty}$ is bounded and we prove that $\rho$ is far from zero and $\|W\|_{H^s}$ is bounded on any bounded time interval $[0,T]$. Using the characteristics 
$$\chi(0,x)\ =\ x, \qquad \chi_t(t,x)\ =\ u(t,\chi(t,x)), $$
the conservation of the mass \eqref{mc} becomes 
$$  \frac{\ud}{\ud\/t}\, \rho\ +\ u_x\/ \rho\  =\ 0, \qquad \implies \qquad   \rho_0(x)\, \mathrm{e}^{-t\, \|u_x\|_{L^\infty}}\ \leqslant\ \rho(t,x)\ \leqslant\ \rho_0(x)\, \mathrm{e}^{t\, \|u_x\|_{L^\infty}},$$
which implies that $\rho$ is bounded and far from zero. The conservation of the energy \eqref{energycons} with the Sobolev embedding $H^1 \hookrightarrow L^\infty$ imply that $\|W\|_{L^\infty}$ is bounded.

Now, we will use that $\rho$ is far from zero and the boundness of $\|W\|_{W^{1, \infty}}$ to prove that $\|W\|_{H^s}$ is also bounded. 
As in the proof of Theorem \ref{energyestimate}, let 
\begin{gather*}
A(W)\ \eqdef\ \begin{pmatrix} \varpi' & 0 \\ 0 & \rho \end{pmatrix}, \qquad   
B(W)\ \eqdef\ \begin{pmatrix} u & \rho \\ \varpi' & u \end{pmatrix}, \qquad
\tilde{E}(W)\ \eqdef\ \left(\Lambda^s\, W,\ A\, \Lambda^s\, W \right),\\
F(W)\ \eqdef\ \begin{pmatrix} 0 \\ -\epsilon\,\mathcal{L}^{-1}_\rho\partial_x \left\{\,\left(\/\rho^2\/\mathscr{A}'\/\right)'u_x^{\,2}\, 
+\left(\/\rho\/\mathscr{V}''\//\/\mathscr{A}'\right)'\,{\mathscr{A}'}^{\/2}\, \rho_x^{\,2}\,\right\} \end{pmatrix},
\end{gather*}
the system \eqref{mc}, \eqref{momG2} then becomes
\begin{equation*}
W_t\ +\ B(W)\, W_x\ =\ F(W).
\end{equation*}
As in \eqref{energyterms}, we have
\begin{align}\label{energyterms2}
\tilde{E}(W)_t\ 
=&\ -2\, \left( [\Lambda^s, B]\, W_x,\, A\, \Lambda^s W \right)\ -\ 2\, \left(B\, \Lambda^s  W_x,\, A\, \Lambda^s W \right) \nonumber \\
&\ - 2\, \left( \Lambda^s F,\, A\, \Lambda^s W \right)\ +\ \left( \Lambda^s W,\, A_t\, \Lambda^s W \right) \nonumber \\
=&:\ I\ +\ II\ +\ III\ +\ IV.
\end{align}
Note that {\small
$$-\half I\, 
=\, \left( [\Lambda^s, u]\, \tilde{\rho}_x, \varpi'\, \tilde{\rho} \right)\,
+\, \left( [\Lambda^s, \rho-\bar{\rho}]\, u_x, \varpi'\, \tilde{\rho} \right)\,
+\, \left( [\Lambda^s, \varpi'(\rho)-\varpi(\bar{\rho})]\, \tilde{\rho}_x, \rho\, u \right)\,
+\, \left( [\Lambda^s, u]\, u_x, \rho\, u \right).$$ }
Using \eqref{Commutator} and \eqref{Composition2} we have 
$$\left| \left( [\Lambda^s, u]\, \tilde{\rho}_x, \varpi'\, \tilde{\rho} \right) \right|\ \lesssim\ \left( \|\tilde{\rho}_x\|_{H^{s-1}}\ +\ \|u\|_{H^s} \right)\, \|\tilde{\rho}\|_{H^s},$$
where the multiplicative constant depend on $\|W\|_{W^{1,\infty}}$. Doing the same for all the terms we obtain that 
\begin{equation}\label{I2}
|I|\ \lesssim\ \|W\|_{H^s}^2\ \lesssim\ \tilde{E}(W).
\end{equation} 
As in \eqref{II}, we obtain 
\begin{equation}\label{II2}
|II|\ \lesssim\ \tilde{E}(W).
\end{equation}
To estimate $III$, we use \eqref{estimate1} and \eqref{estimate4} to obtain that 
\begin{equation}\label{F}
\|F\|_{H^s}\ \lesssim\ \|\psi\|_{H^{s-1}}\ +\ \|\tilde{\rho}\|_{H^{s-1}}\, \left(\|\psi\|_{L^\infty}\, +\,  \|\psi\|_{L^1} \right),
\end{equation}
where $\psi= \left(\/\rho^2\/\mathscr{A}'\/\right)'u_x^{\,2}\, 
+\left(\/\rho\/\mathscr{V}''\//\/\mathscr{A}'\right)'\,{\mathscr{A}'}^{\/2}\, \rho_x^{\,2}$. 
Using \eqref{Algebra}, \eqref{Composition2} one obtain that $ \|\psi\|_{H^{s-1}} \lesssim \|W\|_{H^s}$. 
Due to the conservation of the energy \eqref{energycons}, the quantity $\|W\|_{H^1}$ is bounded, then $\|\psi\|_{L^1}$ is also bounded.
Using that $\|W\|_{W^{1,\infty}}$ is bounded, we obtain that $\|\psi\|_{L^\infty}$ is also bounded. The inequality \eqref{F} then becomes
\begin{equation*}
\|F\|_{H^s}\ \lesssim\ \|W\|_{H^s},
\end{equation*}
which implies that 
\begin{equation}\label{III2}
|III|\ \lesssim\ \|W\|_{H^s}^2\ \lesssim\ \tilde{E}(W).
\end{equation} 
The conservation of the mass \eqref{mc} implies that $\|\rho_t\|_{L^\infty}=\|(\rho u)_x\|_{L^\infty}$ which is bounded. Then 
\begin{equation}\label{IV2}
|IV|\ \lesssim\ \|W\|_{H^s}^2\ \lesssim\ \tilde{E}(W).
\end{equation} 
The equations \eqref{I2}, \eqref{II2}, \eqref{III2}, \eqref{IV2} and \eqref{energyterms2} imply that 
\begin{equation*}
\tilde{E}(W)_t\ \lesssim\ \tilde{E}(W).
\end{equation*}
Gronwall's lemma implies that $\tilde{E}(W)$ does not blow-up in finite time. This ends the proof of the blow-up criteria \eqref{bucrE}.
\qed

\section{A generalised two-component Hunter--Saxton system}\label{secHS}

We have introduced the Sturm--Louville operator
$\mathcal{L}_\rho\/=\/\rho\,-\,2\/\epsilon\/\partial_x\/\rho\/\mathscr{A}'\/\partial_x$
and its inverse $\mathcal{G}_\rho\/=\left[\/1\/-\/2\/\epsilon\/\rho^{-1}\/\partial_x\/
\rho\/\mathscr{A}'\/\partial_x\/\right]^{-1}\/\rho^{-1}$. 
At high frequencies, the operator $\partial_x\mathcal{G}_\rho\partial_x$ behaves like
\begin{equation}\label{Ghinf}
\partial_x\,\mathcal{G}_\rho\,\partial_x\ \sim\ -\/\half\,\epsilon^{-1}\,(\rho\/\mathscr{A}')^{-1}.
\end{equation}
Thus, differentiating with respect of $x$ the equation (\ref{momG}) and considering the high-frequency 
approximation (\ref{Ghinf}), the rbE equations become the system of equations
\begin{align}
\rho_t\ +\,\left[\,\rho\,u\,\right]_x\ &=\ 0, \label{HSgen1}\\ 
\left[\,u_t\,+\,u\,u_x\,+\,\varpi_x\,\right]_x\ &=\, \left(\,1\,+\,\frac{\rho\,\mathscr{A}''}{2\,\mathscr{A}'}
\,\right)u_x^{\,2}\ +\,\left(\,\frac{\left(\rho\/\mathscr{V}''\right)'}{2\,\rho}\,-\,\frac{\mathscr{V}''\,
\mathscr{A}''}{2\,\mathscr{A}'}\,\right)\rho_x^{\,2}, \label{HSgen2}
\end{align}
that is a two-component generalisation of the Hunter--Saxton (gHS) equation \citep{HunterSaxton1991}. 
Smooth solutions of \eqref{HSgen1}, \eqref{HSgen2} satisfy the energy equation 
\begin{equation}\label{HSene}
\left[\rho\, \mathscr{A}'\, u_x^2\ +\ \mathscr{A}'\, \mathscr{V}''\, \rho_x^2  \right]_t\ +\ 
\left[ \left(\rho\, \mathscr{A}'\, u_x^2\, +\, \mathscr{A}'\, \mathscr{V}''\, \rho_x^2 \right)\! u\ +\ 2\, \rho\, \mathscr{A}'\, \mathscr{V}''\, \rho_x\, u_x \right]_x\ =\ 0.
\end{equation} 
There are several generalisations of the Hunter--Saxton equation in the literature, 
including two-component generalisations. The generalisation (\ref{HSgen1})--(\ref{HSgen2}) is 
apparently new and 
it deserves to be studied since it is a simpler system than rbE, being somehow an asymptotic 
approximation. 

It should be noted that equation (\ref{HSgen2}) corresponds to $\mathscr{R}=0$, as easily 
seen considering (\ref{defR}).  
From a physical viewpoint, the Euler equations describe the ``outer'' part  of a shock, while the 
gHS equations describe its ``inner'' structure; the rbE equations being an unification of these 
two (outer and inner) systems. 


Integrating \eqref{HSgen2} with respect to $x$, we obtain 
\begin{align}
\rho_t\ +\,\left[\,\rho\,u\,\right]_x\ &=\ 0, \label{HSgen3}\\ 
u_t\,+\,u\,u_x\,+\,\varpi_x\ 
&=\, \partial_x^{-1}\,  \left\{ \left(\,1\,+\,\frac{\rho\,\mathscr{A}''}{2\,\mathscr{A}'}
\,\right)u_x^{\,2}\ +\,\left(\,\frac{\left(\rho\/\mathscr{V}''\right)'}{2\,\rho}\,-\,\frac{\mathscr{V}''\,
\mathscr{A}''}{2\,\mathscr{A}'}\,\right)\rho_x^{\,2} \right\}\ +\ g(t), \label{HSgen4}
\end{align}
where $(\partial_x^{-1}f)(x) \eqdef \int_0^xf(y)\, \ud y$ and $g(t)=u_t(t,0)+u(t,0)u_x(t,0)+\varpi'(\rho(t,0)) \rho_x(t,0)$.

In the case $\varpi' \equiv 0$, the proof of local well-posedness of \eqref{HSgen3}, \eqref{HSgen4} can be done by using Kato's theorem \cite{kato1975quasi} as in \cite{liu2012,
liu2014,
Moon2012,
Yin2004}. 
Following the proof of Theorem \ref{existencerEuler} and using the inequality 
\begin{equation*}
\|\partial_x^{-1}f\|_{H^{s+1}([0,1])}\ \lesssim\ \|f\|_{H^s([0,1])} \quad \forall s\, \geqslant\, 0,
\end{equation*}
one can prove the following theorem 
\begin{thm}\label{existenceHS}
Let $\tilde{m} \geqslant s \geqslant 2$, $P, \mathscr{A} \in {C}^{\tilde{m}+4}(]0,+\infty[)$ such that $P'(\rho)>0$, $\mathscr{A}'(\rho)>0$ for $\rho>0$. 
Let also $W_0 \in H^s([0,1])$ be a periodic initial data satisfying $\inf_{x \in [0,1]} \rho_0(x) > \rho^{*}$ and $g \in C([0,+\infty[) $, then there exist $T>0$ and a unique periodic solution $W \in {C}([0,T], H^s) \cap {C}^1([0,T], H^{s-1})$ of \eqref{HSgen3}, \eqref{HSgen4} satisfying the non-emptiness condition $ \inf_{x \in [0,1]}\, \rho(t,x)\ >\ 0$ and the conservation of the energy
\begin{equation*}
\frac{\ud}{\ud\/t} \int_0^1 \left( \rho\,\mathscr{A}'\,u_x^{\,2}\ +\ \mathscr{A}'\,\mathscr{V}''\,
\rho_x^{\,2}\,\right)\, \ud x\ =\ 0.
\end{equation*}
Moreover, the maximal existence time $T_{max}< +\infty$, then 
\begin{equation*}
\lim_{t \to T_{max}} \|W_x\|_{L^\infty}\ =\ +\infty.
\end{equation*}
\end{thm}
\begin{rem}
The system \eqref{HSgen3}, \eqref{HSgen4} do not change if $\mathscr{A}$ is replaced by $-\mathscr{A}$.
Then, the result of Theorem \ref{existenceHS} holds in the case $\mathscr{A}'(\rho)<0$. 
\end{rem}

\section{Remarks on a special regularision}\label{secsr}

As proved in \cite{liu2019well}, the solutions given by Theorem \ref{existencerEuler} do not hold for all time in general.
An inspiring way to obtain global (in time) weak solutions, is to use an equivalent semi-linear system of ordinary differential equations as in 
\citep{BressanConstantin2007a,
BressanConstantin2007b,
Grunert2012Global,
Wang2010Global}.
In this case, the lemma \ref{Inverseesitimates} is not enough, and an explicit formula of the operator $\mathcal{L}_\rho^{-1}$ is needed. 

At this stage, the regularising factor $\mathscr{A}$ can be chosen freely, provided that 
$\mathscr{A}'>0$. Here, we investigate further the special choice 
\begin{equation}\label{defApart}
\mathscr{A}(\rho)\ =\ -A\,\bar{\rho}\,/\,\rho,
\end{equation}
where $A>0$ is a constant. 
We show in this section that with this special choice of $\mathscr{A}$, a formula of 
$\mathcal{L}_\rho^{-1}$ can be obtained, and the rbE system can be simplified.
With the special choice \eqref{defApart}, the Sturm--Liouville operator $\mathcal{L}_\rho$ becomes
\begin{equation*}
\mathcal{L}_\rho\ =\ \rho\ -\ 2\,\epsilon\,A\,\bar{\rho}\,\partial_x\,\rho^{-1}\,\partial_x\ 
=\ \rho\left[\,1\,-\,2\,\epsilon\,A\,\bar{\rho}\,\rho^{-1}\,\partial_x\,\rho^{-1}\,
\partial_x\,\right],
\end{equation*}
so its inverse $\mathcal{G}_\rho=\mathcal{L}_\rho^{-1}$ is
\begin{equation*}
\mathcal{G}_\rho\ =\ \left[\,1\,-\,2\,\epsilon\,A\,\bar{\rho}\,\rho^{-1}\,\partial_x\,\rho^{-1}\,
\partial_x\,\right]^{-1}\,\rho^{-1}.
\end{equation*}
Similarly, the operator $\mathcal{J}_\rho$ becomes
\begin{align*}
\mathcal{J}_\rho\ =\ \partial_x^{-1}\,\rho\,\mathcal{G}_\rho\,\partial_x\ =\ \left[\,1\,-\,
2\,\epsilon\,A\,\bar{\rho}\,\rho^{-1}\,\partial_x\,\rho^{-1}\,\partial_x\,\right]^{-1}.
\end{align*}

This special choice for $\mathscr{A}$ suggests the change of independent variables 
$(t,x)\mapsto(\tau,\xi)$ with
\begin{equation*}
\tau\ \eqdef\ t, \qquad 
\xi\ \eqdef\ \int\/\rho(t,x)\,\ud\/x,
\end{equation*}
i.e., $\xi$ is a density potential (defined modulo an arbitrary function of $t$). After 
one spacial integration, the equation (\ref{dLepsdphi}) for the mass conservation yields
\begin{equation*}
\xi_t\ +\ u\,\xi_x\ =\ K(t)\ \equiv\ 0,
\end{equation*}
$K(t)$ being an arbitrary function of $t$ (an integration `constant') that can be set to 
zero without loss of generality, thus providing a gauge condition for $\xi$ (i.e., $\xi$ 
is no longer defined modulo an arbitrary function of $t$).
Thus, with this change of independent variables, the differentiation operators become
\begin{gather*}
\partial_x\ \mapsto\ \rho\,\partial_\xi, \qquad 
\mathcal{J}_\rho\ \mapsto\ \left[\,1\,-\,2\,\epsilon\,A\,\bar{\rho}\,\partial_\xi^{\,2}\,\right]^{-1}, \qquad
\mathcal{G}_\rho\,\partial_x\ \mapsto\ \mathcal{J}_\rho\,\partial_\xi, \\
\partial_t\ \mapsto\ \partial_\tau\ -\ \rho\,u\,\partial_\xi, \qquad
\partial_t\ +\ u\,\partial_x\ \mapsto\ \partial_\tau, 
\end{gather*}
and the regularising term, together with (\ref{defApart}), becomes
\begin{equation}\label{Rpart}
\mathscr{R}\ =\ A\,\bar{\rho}\,\mathfrak{J}\ast\left\{\left(\/\rho\/\mathscr{V}'''\/+
\/3\/\mathscr{V}''\right)\rho_\xi^{\,2}\,\right\}, \qquad
\mathfrak{J}(\xi)\ \eqdef\ \frac{1}{2\,\sqrt{\/2\/\epsilon\/A\/\bar{\rho}\/}}
\,\exp\!\left(\frac{-\/|\xi|}{\sqrt{\/2\/\epsilon\/A\/\bar{\rho}\/}}\right),
\end{equation}
where an asterix denotes a convolution product, i.e., $\mathfrak{J}(\xi)\ast f(\xi)\eqdef
\int_{-\infty}^{\infty}\mathfrak{J}(\xi-\tilde{\rho})\/f(\tilde{\rho})\,\ud\/\tilde{\rho}=\int_{-\infty}^{\infty}
\mathfrak{J}(\tilde{\rho})\/f(\xi-\tilde{\rho})\,\ud\/\tilde{\rho}$ for any function $f$. Note that $\mathfrak{J}$ 
is the pseudo-differential operator $\mathcal{J}_\rho$ rewritten as an integral operator, because  
it is more convenient when applied to weakly regular functions.

With $(\tau,\xi)$ as independent variables, the mass and momentum equations, 
respectively (\ref{dLepsdphi}) and (\ref{eqmomnc}), become 
\begin{equation} \label{sysparru}
\rho_\tau\ +\ \rho^2\,u_\xi\ =\ 0, \qquad 
u_\tau\ +\,\left[\,\rho\,\mathscr{V}'\,-\,\mathscr{V}\,+\,\epsilon\,\mathscr{R}\,\right]_\xi\ =\ 0, 
\end{equation}
with $\mathscr{R}$ given by (\ref{Rpart}). Denoting $\upsilon\eqdef1/\rho$ the specific volume, 
the system (\ref{sysparru}) becomes
\begin{equation*}
\upsilon_\tau\ =\ u_\xi, \qquad 
u_\tau\ =\,\left[\,\frac{\ud\,(\upsilon\/\mathscr{V})}{\ud\/\upsilon}\,+\,\epsilon\,
A\,\bar{\rho}\,\mathfrak{J}\ast\!\left\{\frac{\ud^3\,(\upsilon\/\mathscr{V})}{\ud\/\upsilon^3}\,
\upsilon_\xi^{\,2}\,\right\}\,\right]_\xi. 
\end{equation*}
Eliminating $u$ between these two relations, we obtain
\begin{equation*}
\upsilon_{\tau\tau}\ -\,\left[\,\frac{\ud\,(\upsilon\/\mathscr{V})}{\ud\/\upsilon}\,
\right]_{\xi\xi}\ =\ \epsilon\,A\,\bar{\rho}\,\partial_\xi^{\,2}\,\mathfrak{J}\ast\!\left\{
\frac{\ud^3\,(\upsilon\/\mathscr{V})}{\ud\/\upsilon^3}\,\upsilon_\xi^{\,2}\,\right\}. 
\end{equation*}
At high frequencies, this partial differential equation is approximately
\begin{equation*} 
\upsilon_{\tau\tau}\ -\,\left[\,\frac{\ud\,(\upsilon\/\mathscr{V})}{\ud\/\upsilon}\,
\right]_{\xi\xi}\ \approx\ -\,
\frac{\ud^3\,(\upsilon\/\mathscr{V})}{\ud\/\upsilon^3}\,\frac{\upsilon_\xi^{\,2}}{2}, 
\end{equation*} 
that can be rewritten
\begin{equation} \label{PDEuapp}
\upsilon_{\tau\tau}\ -\ \frac{\ud^2\,(\upsilon\/\mathscr{V})}{\ud\/\upsilon^2}\,
\upsilon_{\xi\xi}\ =\ 
\frac{\ud^3\,(\upsilon\/\mathscr{V})}{\ud\/\upsilon^3}\,\frac{\upsilon_\xi^{\,2}}{2}, 
\end{equation} 
that is a proper hyperbolic partial differential equation if
\begin{equation*}
\frac{\ud^2\,(\upsilon\/\mathscr{V})}{\ud\/\upsilon^2}\ >\ 0.
\end{equation*}
Introducing the velocity $c(\upsilon)\/\eqdef\/\sqrt{\ud^2\/(\upsilon\/\mathscr{V})/\ud\/\upsilon^2}$, 
the equation (\ref{PDEuapp}) is rewritten
\begin{equation*}
\upsilon_{\tau\tau}\ -\ c(\upsilon)^2\,\upsilon_{\xi\xi}\ =\ 
c(\upsilon)\,c'(\upsilon)\,\upsilon_\xi^{\,2}, 
\end{equation*} 
an equation appearing in the theory of liquid crystals, for which smooth solution break down 
in finite time \citep{GlasseyEtAl1997}.

\section{Conclusion and perspectives}

In this paper, we have introduced the regularised barotropic Euler system \eqref{rE2}, inspired by 
the Hamiltonian regularisation of the shallow water (Saint-Venant) system with a constant depth 
introduced in \citep{ClamondDutykh2018a}. The latter work is generalised in two ways: (i) considering 
a more general equation (i.e., barotropic Euler); (ii) introducing a family of regularisations (involving 
an arbitrary function $\mathscr{A}(\rho)$). 

For this system --- and also for the periodic generalised two-component Hunter--Saxton system \eqref{HS} 
--- we prove the local (in time) well-posedness in $H^s$ for $s \geqslant 2$ and a blow-up criteria.
As proven by \citet{liu2019well} for the regularised Saint-Venant equations, those solutions do not exist 
for all time, in general. 

An interesting question that remains open is: {\em Due to the energy equations \eqref{eqene} and \eqref{HSene}, do global weak solutions exist in $H^1$ (or in $\dot{H}^1$ for 
\eqref{HS})?} Two possibilities, that have been used for the Camassa--Holm equation, may also work 
for the systems introduced in the present paper, i.e., using a vanishing viscosity \citep{Guan2011Global} 
or using a semi-linear equivalent system 
\citep{BressanConstantin2007a,BressanConstantin2007b,Grunert2012Global,Wang2010Global}.
Another interesting problem is the study of the limiting cases $\varepsilon \to 0$ and $\varepsilon 
\to \infty$ as in \citep{guelmame2020hamiltonian,guelmame2020global}.

The generalised Hunter--Saxton (gHS) \eqref{HS} system of equations is to rbE what the Hunter--Saxton equation 
\citep{HunterSaxton1991} is to the dispersionless Camassa--Holm equation \citep{CamassaHolm1993}, 
i.e., a ``high frequency limit''. Since the Hunter--Saxton equation is integrable \citep{HunterZheng1994}, 
it is of interest to check if this property is shared with the gHS.

\appendix

\section{Isentropic flows}\label{appisent}

Isentropic motions obey the equation of state 
\begin{equation}\label{eqstatisen}
\rho\,/\,\bar{\rho}\ =\,\left(P\,/\bar{P}\right)^{1/\gamma}, \qquad
P\,/\bar{P}\ =\,\left(\rho\,/\,\bar{\rho}\right)^{\gamma},
\end{equation}
where $\bar{\rho}$ and $\bar{P}$ are positive constants characterising the fluid physical properties at the 
rest state, and $\gamma\eqdef C_\text{p}/C_\text{v}$ is the (constant) ratio of  specific heats 
$C_\text{p}$ and $C_\text{v}$. It should be noted that the constitutive relation (\ref{eqstatisen}) 
gauges the pressure field, so zero pressure level can no longer be chosen arbitrarily  without loss 
of generality. 
For these isentropic motions, we have if $\gamma\neq1$ (taking $\varpi(\bar{P})=0$)
\begin{align}\label{varpiisen}
\varpi\ &=\, \int\left(\frac{\bar{P}}{P}\right)^{\!{1\over\gamma}}\,\frac{\ud\/P}{\bar{\rho}}\ =\ 
\bar{\varpi}\,\frac{(P/\bar{P})^{\frac{\gamma-1}{\gamma}}\,-\,1}{\gamma\,-\,1}\ =\ \bar{\varpi}\,
\frac{(\rho/\bar{\rho})^{\gamma-1}\,-\,1}{\gamma\,-\,1}, \\
\frac{\mathscr{V}}{\bar{P}}\ &=\ \frac{\gamma}{\gamma\,-\,1}\left[\,\frac{1}{\gamma}\left(\frac{\rho}
{\bar{\rho}}\right)^{\!\gamma}\,-\,\frac{\rho}{\bar{\rho}}\,\right]\,=\ \frac{\gamma}{\gamma\,-\,1}
\left[\,\frac{1}{\gamma}\,\frac{P}{\bar{P}}\,-\,\frac{\rho}{\bar{\rho}}\,\right], \nonumber
\end{align}
where $\bar{\varpi}\eqdef\gamma\/\bar{P}\//\/\bar{\rho}$, thence
\begin{equation*}
P/\bar{P}\ =\,\left(\,c_s^{\,2}/\bar{\varpi}\,\right)^{\gamma/(\gamma-1)}, \qquad
c_s^{\,2}/\bar{\varpi}\ =\ 1\ +\ (\gamma-1)\,(\varpi/\bar{\varpi}).
\end{equation*}
In the limiting case $\gamma\to1$ (isothermal motions), these relations become
\begin{equation*} 
c_s^{\,2}\, =\, \bar{\varpi}\, =\, \frac{\bar{P}}{\bar{\rho}}, \quad \frac{\rho}{\bar{\rho}}\, =\, 
\frac{P}{\bar{P}}\, =\, \exp\!\left(\frac{\varpi}{\bar{\varpi}}\right), \quad
\frac{\mathscr{V}}{\bar{P}}\, =\, \frac{\rho}{\bar{\rho}}\log\!\left|\frac{\rho}{\bar{\rho}}\right|
\,-\,\frac{\rho}{\bar{\rho}},
\end{equation*}
so the speed of sound is constant while the density is not. The special case $\gamma=1$ is 
relevant for applications in oceanography because for seawater (at salinity $35\,\mathsf{g/kg}$ 
and atmospheric pressure) $\gamma\approx1.0004$ at $0^\circ\mathsf{C}$ and $\gamma\approx1.0106$ 
at $20^\circ\mathsf{C}$ \citep{KayeLaby1995}.


\end{document}